\documentclass[reqno]{amsart}

\pdfoutput = 1
 
\usepackage[english]{babel}
\usepackage{amsmath,amssymb}
\usepackage{amsfonts}
\usepackage{bm}

\usepackage{amsthm}
\usepackage{array,dcolumn}
\usepackage{graphicx}
\usepackage{epstopdf}

\usepackage{mathtools}
\usepackage[tmargin=3.3cm,bmargin=3.3cm,rmargin=3cm,lmargin=3cm]{geometry}

\usepackage[T1]{fontenc}
\usepackage{lmodern}

\usepackage{caption}
\usepackage[titletoc]{appendix}

\usepackage{subfigure} 

\usepackage{tikz}
\usetikzlibrary{arrows,calc,decorations.pathmorphing}
\usepackage{pgfplots}
\pgfplotsset{filter discard warning=false,compat=newest}

\tikzset{every picture/.style={line width=2pt}}

\usepackage[T1]{fontenc}
\usepackage[utf8]{inputenc}
\usepackage[activate={true,nocompatibility},spacing,kerning]{microtype}
\usepackage[hyperpageref]{backref}
\usepackage[pdftex,
colorlinks,
hyperfootnotes=false,
pdffitwindow=true,
plainpages=false,
pdfpagelabels=true,
pdfpagemode=UseOutlines,
pdfpagelayout=SinglePage,
pdftitle={A locally field-aligned discontinuous Galerkin method for anisotropic wave equations},
pdfauthor={Benedict Dingfelder and Florian Hindenlang,
Max-Planck-Institut für Plasmaphysik, Technische Universität München},
pagebackref,
hyperindex,
filecolor=purple,
urlcolor=purple,
citecolor=blue,
linkcolor=blue]{hyperref}

\newcommand{\textred}[1]{\textcolor{black}{#1}}
\newcommand{\textblue}[1]{\textcolor{black}{#1}}

\renewcommand{\leq}{\leqslant}
\renewcommand{\geq}{\geqslant}

\newcommand{\C}{\mathbb{C}}

\DeclareMathOperator*{\argmax}{argmax}

\newcommand{\Z}{{\mathbb Z}}

\newcommand{\dS}{\ {\mathrm dS}}
\newcommand{\dV}{\ {\mathrm dV}}

\newcommand{\fluxAvg}[1]{ \{\!\!\{ #1 \}\!\!\} }
\newcommand{\fluxJump}[1]{ [\! [ #1 ]\! ] }

\newcommand{\abs}[1]{\left| #1 \right| }
\newcommand{\dof}{\operatorname{DoF}}

\newcommand{\vb}{{\bm{b}}}

\newcommand{\vn}{{\bm{n}}}

\newcommand{\vv}{{\bm{v}}}

\newcommand{\vPhi}{{\bm{\Phi}}}

\newcommand{\vB}{{\bm{B}}}

\newcommand{\vx}{{\bm{x}}}

\newcommand{\Proof}[1]{	\ \\	\textit{Proof:}\\	#1 \\ \parbox[t][2mm]{\textwidth}{ \vspace{-7mm} \flushright $\blacksquare$ }  }

\newtheorem{Theorem}{Theorem}[section]
\newtheorem{Lemma}[Theorem]{Lemma}
\begin{document}

\renewcommand{\figurename}{\small {\sc Figure\@}}
\renewcommand{\tablename}{\small {\sc Table\@}}

%%=====================
%%svjour
%%=====================
%
\title[\tiny{A locally field-aligned discontinuous Galerkin method for anisotropic wave equations}]{A locally field-aligned discontinuous Galerkin method for anisotropic wave equations}
%
%\author{Benedict Dingfelder \and Florian J. Hindenlang}
%
%
%\institute{Benedict Dingfelder (\email{benedict.dingfelder@ipp.mpg.de})
%           \and Florian J. Hindenlang (\email{florian.hindenlang@ipp.mpg.de}) \at Max Planck Institute for Plasma Physics, Boltzmannstra{\ss}e 2, D-85748 Garching, Germany}
% 
%\maketitle

%=====================
%amsart 
%=====================

 \author{Benedict Dingfelder}
 \address{\scriptsize{Max-Planck-Institut für Plasmaphysik, 85748~Garching, Germany}}
 \email{benedict.dingfelder@ipp.mpg.de}
 \author{Florian J. Hindenlang}
 \address{\scriptsize{Max-Planck-Institut für Plasmaphysik, 85748~Garching, Germany}}
 \email{florian.hindenlang@ipp.mpg.de}

%%\subjclass[2010]{65E05, 65D25; 68R10, 05C38}

\begin{abstract}
	In magnetized plasmas of fusion devices the strong magnetic field leads to highly anisotropic physics where solution scales along field lines are much larger than perpendicular to it.
	Hence, regarding both accuracy and efficiency, a numerical method should allow to address parallel and perpendicular resolutions independently. In this work, we consider the eigenvalue problem of a two-dimensional anisotropic wave equation with variable coefficients which is a simplified model of linearized ideal magnetohydrodynamics.\\
	For this, we propose to use a mesh that is aligned with the magnetic field and choose to discretize the problem with a discontinuous Galerkin method which naturally allows for non-conforming interfaces.\\
	First, we analyze the eigenvalue spectrum of a constant coefficient anisotropic wave equation, and demonstrate that this approach improves the accuracy by up to seven orders of magnitude, if compared to a non-aligned method with the same number of degrees of freedom. In particular, the results improve for eigenfunctions with high mode numbers.\\
	We also apply the method to compute the eigenvalue spectrum of the associated anisotropic wave equation with variable coefficients of flux surfaces of a Stellarator configuration. We benchmark the results against a spectral code.
%	This provides the foundation to further extend the framework to full MHD problems and three-dimensional locally field-aligned hexahedral elements.

%	\TODO{not up to date}
%	In magnetized plasmas of fusion devices the strong magnetic field leads to highly anisotropic physics. 
%	If only diffusion processes are considered, the diffusion along the magnetic field is dominating. In the limit of vanishing perpendicular diffusion, we obtain the anisotropic diffusion equation with a semidefinite diffusion tensor. We consider the associated eigenvalue problem given by
%	\begin{equation*}
%		-\nabla \cdot \left( \vb \vb \cdot \nabla \phi \right) = \omega^2 \phi \qquad \text{ in } \Omega \subset \R^2
%	\end{equation*}
%	for the two-dimensional fully periodic domain $\Omega$ and direction of the magnetic field $\vb$. \textbf{TODO: This eigenvalue problem is difficult to solve as the differential operator is non-coercive.} We propose a discontinuous Galerkin (DG) method on a non-conforming mesh with locally aligned cells which improves the numerical accuracy by multiple orders of magnitude in comparison to non-aligned methods with the same computational complexity.
\end{abstract}

\keywords{Discontinuous Galerkin, field-alignment, non-conforming, anisotropic wave, non-coercive, eigenvalue problem, ideal MHD, plasma physics, Stellarator, DOI: \href{https://doi.org/10.1016/j.jcp.2020.109273}{https://doi.org/10.1016/j.jcp.2020.109273}}

%\category{G.1.2}{Numerical Analysis}{Approximation, quadrature}
%
%\terms{Inverse Laplace transform, Talbot's method, trapezoidal rule, midpoint rule}
%

%\tableofcontents

%\begin{bottomstuff}
%The research of BD is supported by ({\sl Benedict please}).
%The research of JACW is supported by the National Research
%Foundation of South Africa.
%
%Author's addresses: B.~Dingfelder, Center for Mathematics,
%Technical University of Munich, Boltzmannstrasse 3,
%85748 Garching bei M\"unchen, Germany ({\tt dingfelder@mytum.de});  J. A. C.~Weideman,
%Department of Mathematical Sciences, Stellenbosch University,
%Stellenbosch 7600, South Africa ({\tt weideman@sun.ac.za}).
%\end{bottomstuff}

\maketitle

\section{Introduction}  
\label{sec:intro}
% 
% \begin{itemize}
% 	\item Description of the problem.
% 	\item Why are commonly used methods problematic for this equation?
% 	\item Physical motivation for reducing the space of modes.
% 	\item motivation for small eigenvalues $\rightarrow$ low frequencies
% 	\item motivation for maximal error over all eigenvalues
% 	
% 	\item Flux surfaces of MHD equilibria
% 	\item Straight-field coordinate system in a flux surface, $\iota$ profile
% 	
% 	\item define what is ADG, FE, etc. (if needed)
% \end{itemize}
% 

In magnetically confined fusion devices, hot plasma is confined by strong magnetic fields. Two main torus-shaped designs are distinguished, the Tokamak with an axisymmetric field and the Stellarator with a fully three-dimensional field. In general, the magnetic field geometry of a confined state can be described by a so-called MHD equilibrium, a non trivial steady-state of the ideal MHD equations where magnetic and pressure forces balance \cite{freidberg}. 
\begin{figure}[!htb]
	\centering
	\begin{tikzpicture}[line width =1pt]
	\draw (-4.5, 1.1) node[inner sep=0] {
		\includegraphics[width=0.52\textwidth]{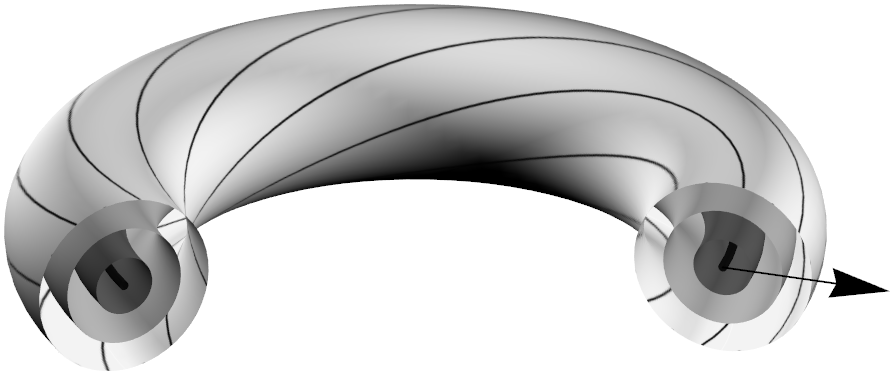}
	};
	\draw (-0.65,0.45) node {$s$};
	\draw[->] (0,0) arc (-45:45:1);	
	\draw[above] (0,1.4) node {$\varphi$};
	\draw[line width=5pt,draw=white,->] (0.1,-0.15) arc (-90:180:0.33);
	\draw[->] (0.1,-0.15) arc (-90:180:0.33);
	\draw[below] (0.1-0.33,-0.1+0.33) node {$\theta$};
	
	\end{tikzpicture}
	\caption{Sketch of nested flux surfaces in a toroidal geometry. Magnetic field lines are shown on the outermost surface.} 
	\label{fig:torus_flux_surfaces}
\end{figure}
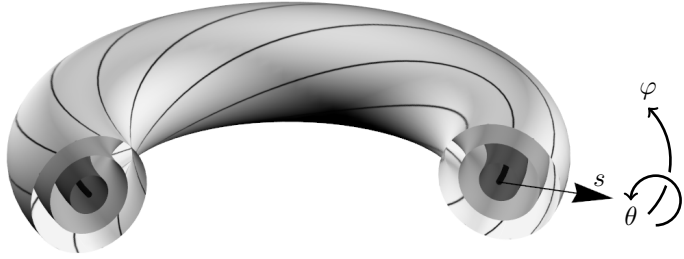\\
In the core of the plasma, the magnetic field lines lie on a set of nested toroidal surfaces of constant pressure, so-called magnetic flux surfaces, as depicted in Figure~\ref{fig:torus_flux_surfaces}. We introduce a radial coordinate $s\in[0,1]$ that labels the flux surfaces from the magnetic axis to the largest flux surface. The flux surface geometry is parametrized by two angles, the poloidal angle $\theta$ and toroidal angle $\varphi$. The magnetic field is tangential to the flux surface, thus it is defined by the flux surface geometry and the contra-variant components $(B^\theta ,B^\varphi)$. It is now possible to define straight-field line coordinates (SFL), or magnetic coordinates, in which the magnetic field lines are straight~\cite{SFLcoordinates}. Hence, the ratio of the contra-variant components is a constant, called the rotational transform $\iota$, and only depends on the flux surface label
\begin{equation}
 \iota(s):=B^\theta(s,\theta,\varphi)/B^\varphi(s,\theta,\varphi) \ . \label{eq:iota_def}
\end{equation} 
The magnetic field geometry of an MHD equilibrium is the starting point for studying the physical behavior of the plasma. For example, finding the associated eigenmode structure of the linearized MHD equations will give insight into the stability of the equilibrium against disturbances and also about resonance modes~\cite{freidberg}. Due to the strong magnetic field, the plasma response is highly anisotropic, such that gradients in magnetic field direction are typically much smaller than in the perpendicular direction.\\
Hence, regarding both accuracy and efficiency, a numerical method should allow to address parallel and perpendicular resolutions independently.\\
One approach is to use finite difference methods and trace the solution along magnetic field lines to approximate parallel derivatives. Amongst others, \cite{fd_beginnings,bruce_fd,fd_overview,grillix,ericGyrokinetic} study such
finite difference approaches for plasma turbulence. In this work, we want to mimic such a strategy with a mesh-based finite element method. Therefore, the element edges of the mesh have to follow the magnetic field. Global alignment is not possible as magnetic field lines do not necessarily close, so we propose to align elements only locally between two consecutive poloidal planes. Therefore, non-conforming interfaces have to be introduced. Even though mortar methods for finite elements exist \cite{mortar_FE}, we choose to construct a high order discontinuous
Galerkin (DG) method which naturally incorporates the treatment of
non-conforming interfaces \cite{mortar_DG}.\\
\textblue{The solution of eigenvalue problems using Finite element methods has been widely studied, an overview is given by Boffi \cite{boffi2010}. Regarding the DG method with conforming and non-conforming interfaces, theoretical bounds and spectrally correct solutions have been obtained for the eigenproblem of the Laplace operator in \cite{antonietti2006} and, in particular, the eigenproblem of the Maxwell equations (curl-curl operator). The Maxwell eigenproblem is first studied in \cite{hesthaven2004_DG_Eigval_maxwell,warburton2006_LDG_eigval_maxwell} for  applications on conforming meshes, showing spectrally correct solutions, if the penalty parameter is chosen sufficiently high. In \cite{buffa2006,buffa2007}, the theoretical analysis is conducted in more detail for non-conforming meshes,  and in \cite{buffa2009}, a spectrally correct mortar-type DG method is introduced for general non-conforming interfaces.
Both the curl-curl operator and the anisotropic wave equation considered here are non-coercive operators. However, as will be shown in Section~\ref{sec:eigenfunctions_gradients}, an additional complexity of the anisotropic wave equation is that arbitrarily small eigenvalues exist, even for eigenvectors with large mode numbers.  }\\
\textblue{Our contribution is the adaptation of the local DG method for the anisotropic wave equation, using a mixed variational form and a locally aligned non-conforming mesh. We show symmetry of the discretized system and the numerical convergence of the eigenvalues for constant and variable coefficients.   
%
%We propose to align the elements only locally, so that non-conforming interfaces have to be introduced. 
%
%Furthermore, discontinuous Galerkin methods offer great freedom in the
%subdivision of the domain as well as the choice of the basis for each
%cell. 
Further, we demonstrate that the proposed DG method on the aligned mesh is able to address the degrees of
freedom for resolving parallel and perpendicular direction individually. 
}\\
%To the best of our knowledge, this is the first time that a mesh-based finite element method with a locally aligned mesh is applied to a highly anisotropic problem.
%
As a model problem, we study a simplified linear MHD model that keeps the anisotropy of the underlying physics and describes the resonance behavior of the plasma on a single flux surface. It is a two-dimensional anisotropic wave equation with variable coefficients that arise from the geometry of the magnetic flux surface and was derived in \cite{mydiss}. For the specific case of constant coefficients, the eigenmodes can be computed analytically, allowing us to show a distinct feature of the spectrum that cannot be seen in simple resonant systems. Low frequencies are not directly linked to low mode numbers, but also high mode numbers can have a low frequency, if the associated eigenfunction is aligned with the magnetic field.\\
%
%For simple resonant systems, low mode numbers can be directly associated to low frequencies. However, a
%A main feature of the spectrum of the anisotropic problem is that also high mode numbers can have a low frequency, if the associated eigenfunction is aligned with the magnetic field. 
In Section~\ref{sec:equation}, the equation and the analytic properties with constant coefficients are discussed. 
In Section~\ref{sec:discretization}, we show how the aligned mesh is constructed and derive the discontinuous Galerkin method for the anisotropic wave equation with variable coefficients. In Section~\ref{sec:results}, we investigate the properties of the proposed method for the case of constant coefficients. Here, the errors to exact eigenvalues can be evaluated, allowing to  quantify the impact of the alignment, the impact of the resolution ratio between parallel and perpendicular direction and also show the mesh convergence of the method. Finally, in Section~\ref{sec:ADM_MHD}, we consider a three-dimensional magnetic field geometry of a Stellarator, with variable coefficients for each flux surface. We demonstrate the convergence of the eigenvalues and benchmark the result against a high resolution spectral code.

\section{Anisotropic wave equation}
\label{sec:equation}
The two-dimensional anisotropic wave equation with variable coefficients  was derived in~\cite{mydiss} from the linearized ideal MHD equations, with drift approximation~\cite{bruce_drift_approximation} and by neglecting  pressure and compressional Alfvén waves.\\
The associated two-dimensional eigenvalue problem reads as
\begin{equation}
\label{eq:eigenvalue_equation}
-\nabla \cdot \left( \vB \left(\vB \cdot \nabla \phi \right)\right) = \omega^2 \alpha  \phi \ ,
\end{equation}
within the periodic domain $\vx \in \Omega = [0,2\pi)^2$ and the magnetic field 
\begin{equation}
 \vB  = \beta(\vx) \vb = \beta(\vx)  \begin{pmatrix}
 b_1 \\ b_2
 \end{pmatrix} \ .
\end{equation} 
As stated in the introduction, using straight field line coordinates $(x,y)= (\theta,\varphi)$ allows to express the magnetic field by a real valued constant vector $(b_1,b_2)=(\iota,1)$ and a scalar-valued periodic function $\beta\left( \vx \right)>0$. The right hand side is allowed to include a scalar-valued periodic variation $\alpha \left( \vx \right)>0$. Both scalar fields $\alpha,\beta$ arise from the MHD equilibrium and the metric terms of the magnetic flux surface geometry, details are given in~\cite{mydiss}. The special case of a flux surface in a periodic cylinder geometry can be described by $\alpha,\beta\equiv 1$. 

\subsection{Analytical properties for constant coefficients}
\label{sec:eigenfunctions_gradients}
To investigate properties of solutions of \eqref{eq:eigenvalue_equation}, we consider analytic eigenfunctions associated to eigenvalues of \eqref{eq:eigenvalue_equation} for constant coefficients $\alpha, \beta \equiv 1$. Defining the norm
\begin{equation}
\| v \|^2_{L_2(\Omega)} \coloneqq \int_{\Omega} v \bar{v} \dV, 
\end{equation} 
\textred{we }\textblue{observe }\textred{that }
	\begin{equation}
\| \vb \cdot \nabla \phi_{m,n}(x,y) \|^2_{L^2_2(\Omega)}  \leq  4 \pi^2  \omega_{m,n}^2.
\end{equation}
\textblue{where } 
\begin{equation}
\label{eq:exact_eigenfunctions}
\phi_{m,n}(x,y) = \exp \left( \mathrm{i}\left( mx + ny\right)\right), \qquad m,n\in\Z ,\ x,y\in \Omega
\end{equation}
\textred{are }\textblue{the }\textred{analytic eigenfunctions corresponding to the eigenvalue $\omega_{m,n}^2 =(b_1 m + b_2n)^2$. This shows that the gradients of eigenfunctions with small eigenvalues are small along $\vb$. Therefore, the eigenfunctions themselves are close to constant in $\vb$-direction. Hence, less resolution in parallel than in perpendicular direction is needed in a discretization.\\
This is illustrated in Figure~\ref{fig:density_plot_-4_5_mode}, where we show the density plot of the real part of $\phi_{4,-5}$ which has the eigenvalue $\omega_{4,-5}=0.11305798$ for the choice of \textred{$\vb = \left( 1.165939761, 1\right)^\top $}. We observe that the function is almost constant in $\vb$-direction which is plotted as the black dashed line.
\begin{figure}[!htb]
	\centering
	\includegraphics[width=0.7\textwidth]{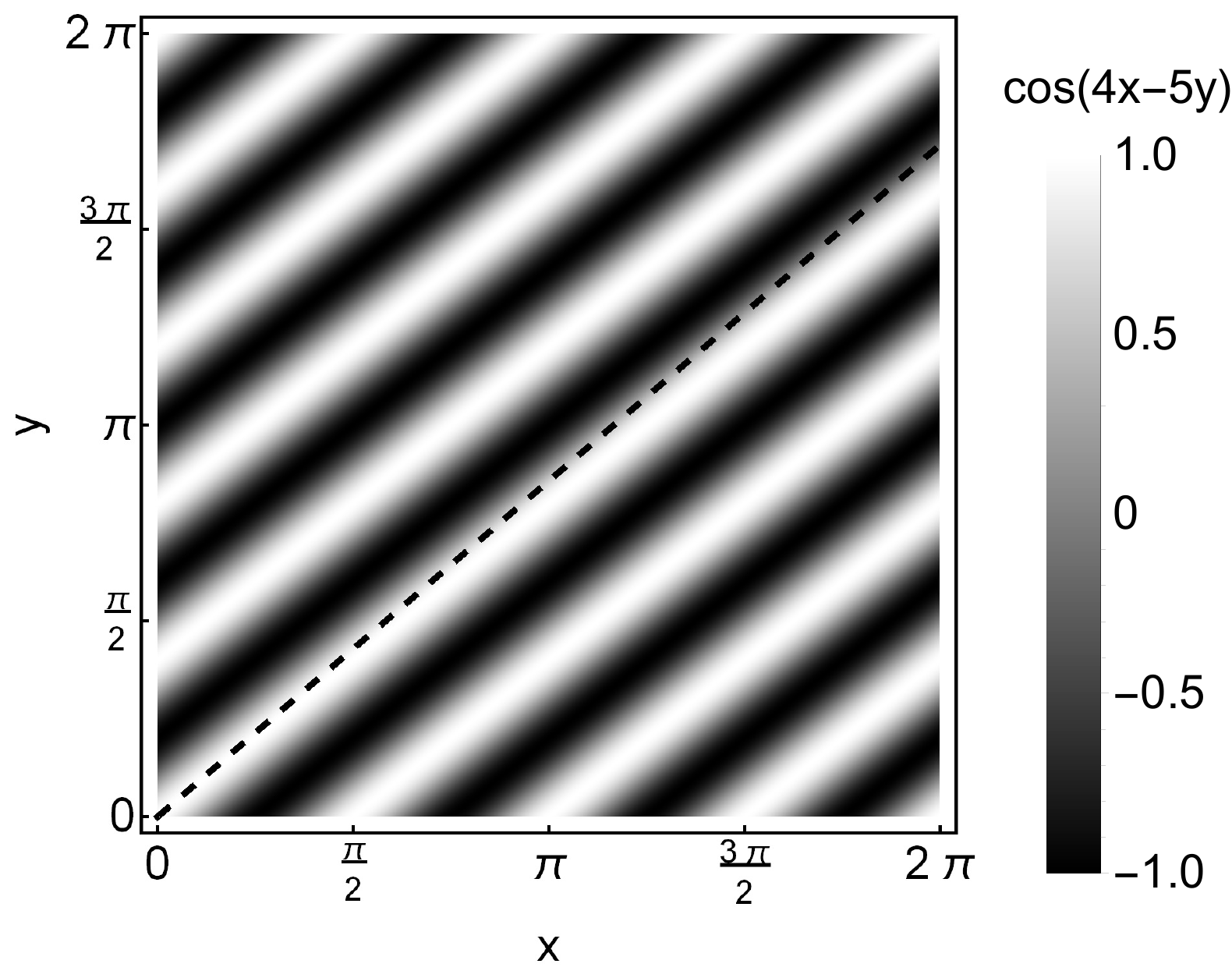}
	\caption{Real part of eigenmode $\phi_{4,-5}$ in $\Omega = [0,2\pi)^2$ with associated eigenvalue $\omega_{4,-5}=0.11305798$ for \textred{$\vb = \left( 1.165939761, 1\right)^\top $}  (black dashed line).} 
	\label{fig:density_plot_-4_5_mode}
\end{figure}\\
}
\textred{We remark that the Fourier modes $\phi_{m,n}$ with mode ratio $m/n$ close to $-b_2/b_1$ are those producing small eigenvalues, as 
\begin{equation}
\omega_{m,n}^2 = (b_1 m + b_2 n)^2 = n^2 (b_1 \frac{m}{n}+b_2)^2 \ .
\end{equation}
This demonstrates that also high mode numbers can be associated with small eigenvalues, if their eigenfunction has little variation in field direction. }

\section{Discretization}
\label{sec:discretization}
%
%In Section~\ref{sec:eigenfunctions_gradients}, we examine the analytical solutions of the eigenvalue problem \eqref{eq:eigenvalue_equation} for $\alpha, \beta \equiv 1$ for determining which features the method needs to resolve. 
In Section~\ref{sec:choice_of_mesh_and_basis}, we \textblue{propose the design of a suited mesh that is aligned with $\vb$. Further, we discuss remarks on the function spaces for discretization of the test and trial functions}. In Section~\ref{sec:mixed_form}, we formulate the discontinuous Galerkin method based on a variational mixed formulation of \eqref{eq:eigenvalue_equation} and derive the discrete eigenvalue problem.

\subsection{Choice of mesh and basis}
\label{sec:choice_of_mesh_and_basis}
\begin{figure}[!htb]
	\centering
	\begin{tikzpicture}[scale = 3]
	\draw[->] (-1/4,0) -- ++(0.75,0);
	\node[right] at (0.5,0) {$x$};
	\draw[->] (-1/4,0) -- ++(0,0.75);
	\node[above] at (-1/4,0.75) {$y$};
	\draw[thick,->] (-1/4,0) -- ++(1/2,1/6);
	\node[below, right] at (1/4,1/6) {$\vb$};	
	\draw[->] (2.5,0.5) -- (3.25,0.5);
	\node at (3.4,0.5) {$\xi$};
	\draw[->] (2.5,0.5) -- (2.5,1.25);
	\node at (2.5,1.4) {$\eta$};
	\draw (0,0.25) -- ++(1,1/3) -- ++(0,1) -- ++(-1,-1/3) -- cycle;
	\node at (0.5-0.25,0.25+1/6+1/2-0.25-1/12) {$K$}; 
	\draw[->] (0.5,0.25+1/6+1/2) -- ++(0.75,1/6+1/12); 
	\draw[->] (0.5,0.25+1/6+1/2) -- ++(0,0.75); 
	\node[right] at (0.5+0.75,0.25+1/6+1/2+1/6+1/12) {$\xi$}; 
	\node[above] at (0.5,0.25+1/6+1/2+0.75) {$\eta$}; 
	\draw (2,0) -- (3,0) -- (3,1) -- (2,1) -- (2,0);
	\node at (1.9,-0.1) {$(-1,-1)$};
	\node at (3.1,-0.1) {$(1,-1)$};
	\node at (3.1,1.1) {$(1,1)$};
	\node at (1.9,1.1) {$(-1,1)$};
	\draw[->] (1.25,0.75) -- (1.75,0.75);
	%	\node at (1.5,0.9) {$\vL^{-1}_K$};
	\draw[<-] (1.25,0.25) -- (1.75,0.25);
	%	\node at (1.5,0.1) {$\vL_K$};
	\end{tikzpicture}\\	
	\caption{\textblue{Reference coordinate system $\left( \xi,\eta \right) $ for an aligned cell $K$ in $\left( x,y\right)$.}}
	\label{fig:CODE_ref_element_map_and_inverse}
\end{figure}
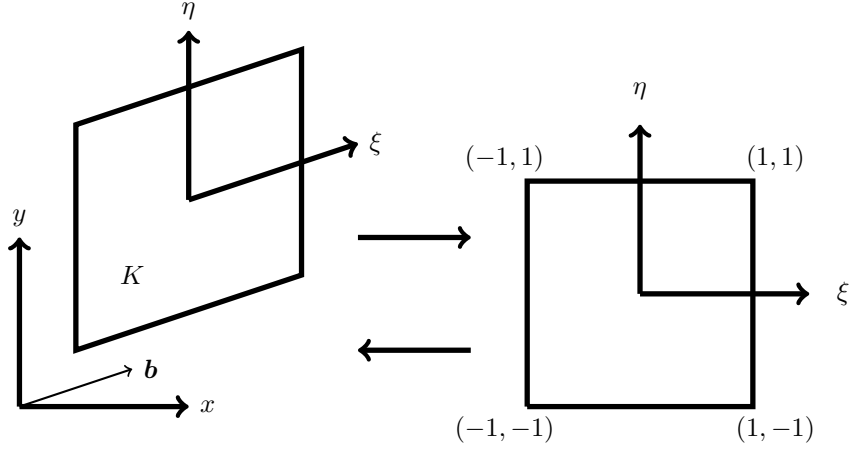

Section~\ref{sec:eigenfunctions_gradients} shows that a suited method should distribute its degrees of freedom and emphasize on the resolution in perpendicular direction. When discretizing the two-dimensional fully periodic domain $\Omega$, we want to design a mesh which aims at uniformly treating all cells and allows us to assign the resolution separately. 
\begin{figure}
	\centering
	\subfigure[cartesian, conforming,  commonly used for fully periodic domain]{			
			\begin{tikzpicture}[scale = 0.35]
			\foreach \k in {1,2,3,4,5,6} {
				\foreach \i in {1,2,3,4,5,6} {
					\draw (\i+0,\k+0) -- (\i+1,\k+0) -- (\i+1,\k+1) -- (\i+0,\k+1) -- (\i+0,\k+0) ;
				}
				
			}
			\draw[white] (-1,0)--(8,0);
			\end{tikzpicture}}
		\hspace{2mm}
	\subfigure[fully aligned, internally conform, periodically non-conforming]{
			\begin{tikzpicture}[scale = 0.35]	
			\foreach \k in {1,2,3,4,5,6} {
				\foreach \i in {1,2,3,4,5,6} {
					\draw (\i+0,\k+\i*1/3-1/3) -- (\i+1,\k+1/3+\i*1/3-1/3) -- (\i+1,\k+4/3+\i*1/3-1/3) -- (\i+0,\k+1+\i*1/3-1/3) -- (\i+0,\k+0+\i*1/3-1/3) ;	
				}
			}	
			\draw[thick,->] (4,1) -- (6,5/3);
			\node at (5.5,0.9) {$\vb$};	
			\draw[white] (-1,0)--(8,0);
			\end{tikzpicture}
	} 
		\hspace{2mm}
	\subfigure[locally aligned, non-conforming, uniform $N_x = N_y$]{
		\begin{tikzpicture}[scale = 0.35]	
			\foreach \k in {1,2,3,4,5,6} {
				\foreach \i in {1,2,3,4,5,6} {
					\draw (\i+0,\k+0) -- (\i+1,\k+1/3) -- (\i+1,\k+4/3) -- (\i+0,\k+1) -- (\i+0,\k+0) ;
				}
			}	
			\draw[white] (-1,0)--(8,0);
		\end{tikzpicture}
		}
		\hspace{2mm}
	\subfigure[locally aligned, non-conforming, generalized $N_x \neq N_y$, here: $N_x = 3, N_y = 6$]{
			\begin{tikzpicture}[scale = 0.35]	
			\foreach \k in {1,2,3,4,5,6} {
				\foreach \i in {1,3,5} {
					\draw (\i+0,\k+0) -- (\i+2,\k+2/3) -- (\i+2,\k+5/3) -- (\i+0,\k+1) -- (\i+0,\k+0) ;
				}
			}
			\draw[white] (-1,0)--(8,0);
			\end{tikzpicture}						
	}
\caption{Different conforming and non-conforming choices for discretizing the domain $\Omega$.}
\label{fig:discretization_meshes}
\end{figure}
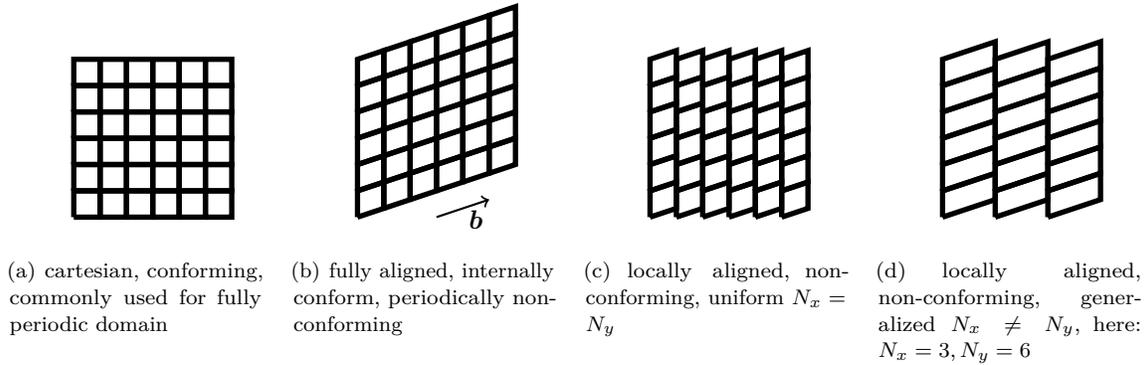\\
\begin{figure}
	\centering
	\subfigure[upper and bottom aligned cells]{
		\begin{tikzpicture}[scale=0.7]	
		\foreach \k in {1,2,3,4} {
			\foreach \i in {1,2,3,4} {
				\draw (\i+0,\k+0) -- ++(1,7/2) -- ++(0,1) -- ++(-1,-7/2) -- cycle ;
			}
		}
		\end{tikzpicture}					
	}
	\subfigure[left and right aligned cells]{
		\begin{tikzpicture}[scale=0.7]	
		\foreach \k in {1,2,3,4} {
			\foreach \i in {1,2,3,4} {
				\draw (\i+0,\k+0) -- ++(1,0) -- ++(2/7,1) -- ++(-1,0) -- cycle ;
			}
		}
		\draw[->] (-1/2,3/2) -- ++(1,3/2+7/4);
		\node[below] at (-1/2,3/2) {$\vb$}; 
		\end{tikzpicture}					
	}
	
	\caption{Aligned meshes  for $ \frac{b_1}{b_2} = \frac{2}{7}$, with large aspect ratio (left) and small aspect ratio (right), for $N_x=N_y=4$. }
	\label{fig:AD_locally_aligned_mesh_strongly_sheared}	
\end{figure}
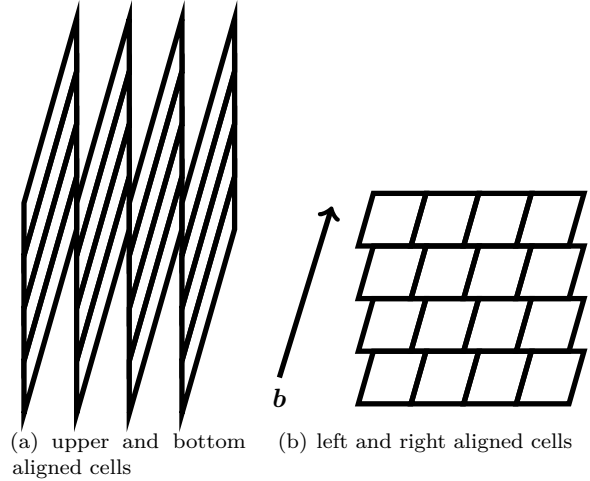
The most intuitive choice is to fully align one dimension of the mesh with $\vb$. This yields non-conforming interfaces at the periodic boundary whenever
\begin{equation}
	\label{eq:periodic_boundary_constraint}
	\frac{b_2}{b_1} N_y \notin \Z
\end{equation}
where $N_x, N_y$ is the number of cells in $x$,$y$-direction respectively. \eqref{eq:periodic_boundary_constraint} particularly holds for all irrational fractions $b_2/b_1$. To generate a more regular distribution of non-conforming interfaces, we choose to align all mesh cells locally instead of globally which introduces non-conforming interfaces for each cell whenever
\begin{equation}
\label{eq:non_conformity_constraint}
\frac{b_2}{b_1} \frac{N_y}{N_x} \notin \Z . 
\end{equation}
Now both, interior and boundary interfaces, are treated equally. Once an aligned mesh is used, we can adapt the resolution in parallel and perpendicular direction by changing $N_x$ and $N_y$ respectively. \\
Figure~\ref{fig:discretization_meshes} shows aligned upper and lower cell interfaces. \textred{ In the case of $\left| b_1 \right|  \ll \left|b_2 \right|$, an alignment of left and right cell interfaces should be chosen instead, as shown in Figure~\ref{fig:AD_locally_aligned_mesh_strongly_sheared}, because the aspect ratio (AR, ratio of longest to shortest side) is very large for the case of aligned upper and lower cell interfaces. The aspect ratio for both choices is proportional to
\begin{equation}
  \textrm{AR}_\textrm{upper/lower} \sim \sqrt{1+\left(\frac{b_2}{b_1}\right)^2}\,,\quad  
  \textrm{AR}_\textrm{left/right} \sim \sqrt{1+\left(\frac{b_1}{b_2}\right)^2}\,.
\end{equation} 
Very large aspect ratios can deteriorate the numerical accuracy of the simulations, as the element Jacobian enters in the condition number of the matrix system.}\\
\textblue { We now consider the discretization of the spaces of test and trial functions for the discontinuous Galerkin method. The basis functions are defined on a reference element $\left( \xi, \eta \right) \in [-1,1]^2$ as a tensor product of two polynomials with degrees $p_\xi,p_\eta$. We choose to align $\xi$ with $\vb$ and $\eta$ with $y$ as depicted in Figure~\ref{fig:CODE_ref_element_map_and_inverse}. This allows us to separate resolutions, using the parameters 
\begin{itemize}
	\item resolution of the mesh $N_x, N_y$ 
	\item  degree of the basis $p_\xi, p_\eta$ 
\end{itemize}
with respective parallel and perpendicular resolution
\begin{equation}
\label{eq:Results_dof_definition}
\dof_\parallel \coloneqq \left( p_\xi +1 \right)  N_x , \qquad  \dof_\perp \coloneqq \left( p_\eta + 1  \right) N_y \ . 
\end{equation}}

\subsection{Mixed form and derivation of the DG method}
\label{sec:mixed_form}

As carried out in \cite[Sections 1.3, 4.2.1]{diPietroErn}, the bilinear form for solving the eigenvalue problem should be consistent, continuous, coercive and symmetric. As the differential operator of \eqref{eq:eigenvalue_equation} is not coercive, this property cannot be transferred to the associated bilinear form. For preserving symmetry, we propose the mixed form
\begin{equation}
\label{eq:symm_mixed_form}
\begin{dcases}
\vB \cdot \nabla \phi &= u \\
-\nabla \cdot \left( \vB u \right) &= \omega^2 \alpha \phi \ . 
\end{dcases}
\end{equation}
Note that the symmetric splitting leads to a scalar-valued equation for the parallel gradient represented by $u$. The associated weak form for test and trial functions $v,u \in V_{U}$, $\psi, \phi \in V_{\Phi}$ and a mesh with cells $K \in \mathcal{K}$ writes
\begin{equation}
\label{eq:AD_mixed_form_var_form_general}
\begin{dcases}
\sum_{K \in \mathcal{K}} \left( \int_K \vB \cdot \nabla \phi v \dV \right)  =  \sum_{K \in \mathcal{K}}\left( \int_K u v \dV \right)  \qquad &\forall v \in V_{U}  \\
\sum_{K \in \mathcal{K}} \left( \int_K -\nabla \cdot \left( \vB u \right) \psi \dV \right) = \omega^2 \sum_{K \in \mathcal{K}} \left( \int_K \alpha \phi \psi \dV \right)  \qquad &\forall \psi \in V_{\Phi}  \ .
\end{dcases}
\end{equation}
Using locally defined test and trial functions $v^K,u^K \in V_{K,U}$, $\psi^K, \phi^K \in V_{K,\Phi}$ in the locally defined function spaces $V_{K,U}, V_{K,\Phi}$ for a cell $K$, \eqref{eq:AD_mixed_form_var_form_general} writes
\begin{equation}
\label{eq:AD_mixed_form_var_form_local}
\begin{dcases}
\int_K \vB \cdot \nabla \phi^K v^K \dV =  \int_K u^K v^K \dV \qquad &\forall v^K \in V_{K,U}  \\
\int_K -\nabla \cdot \left( \vB u^K \right) \psi^K \dV = \omega^2 \int_K \alpha \phi^K \psi^K \dV  \qquad &\forall \psi^K \in V_{K,\Phi} \ . 
\end{dcases}
\end{equation} 
In the following, we omit the spaces for test functions and the superscripts $K$ whenever there is no need for differentiating the cells of definition.
Integration by parts of \eqref{eq:AD_mixed_form_var_form_local} leads to 
\begin{equation}
\begin{dcases}
- \int_K \phi \nabla \cdot \left( \vB  v \right)  \dV  + \sum_{F \in \partial K} \left( \int_{F} \widehat{\phi} v \vB \cdot \vn \dS \right) &=  \int_K u v \dV  \\
\int_K u \vB \cdot \nabla \psi \dV  - \sum_{F \in \partial K} \left( \int_{F} \widehat{u} \psi \vB \cdot \vn \dS \right) & = \omega^2 \int_K \alpha \phi \psi \dV
\end{dcases}
\end{equation}
with numerical fluxes $\widehat{u}$ and $\widehat{\phi}$ that still need to be defined and unit outer normal $\vn$ on $K$. %As the first integrals are not symmetric to each other in the sense of exchanging $u$ with $v$ and $\phi$ with $\psi$, we use integration by parts for the second equation once more. The boundary term for this integration is carried out using the contribution of the function on the main cell, namely $\phi^-$
To recover the symmetry of the volume integrals, the first equation is integrated by parts again using the inner cell boundary value $\phi$, yielding
\begin{equation}
\label{eq:variational_mixed_form}
\begin{dcases}
-\int_K u v \dV  + &\int_K \vB \cdot \nabla \left(  \phi \right) v \dV  + \sum_{F \in \partial K} \left( \int_{F} \left( \widehat{\phi} -\phi \right)  v \vB \cdot \vn \dS \right)  =0 \\
&\int_K u \vB \cdot \nabla \psi \dV - \sum_{F \in \partial K} \left( \int_{F}  \widehat{u} \psi \vB \cdot \vn \dS \right) = \omega^2 \int_K \alpha \phi \psi \dV \ . 
\end{dcases}
\end{equation}
%When discretizing \eqref{eq:variational_mixed_form}, we obtain a linear system with coefficient vectors $U$ and $\Phi$
%\begin{equation}
%\label{eq:mixed_form_system}
%\begin{pmatrix}
%M_1 & M_2 \\
%M_3 & M_4 
%\end{pmatrix}
%\begin{pmatrix}
%U\\
%\Phi
%\end{pmatrix}
%= \omega^2
%\begin{pmatrix}
%0 & 0 \\
%0 & M_5 
%\end{pmatrix}
%\begin{pmatrix}
%U\\
%\Phi
%\end{pmatrix}.
%\end{equation}
%The system can be reduced using $U = -M_1^{-1} M_2 \Phi$ to
%\begin{equation}
% \left( - M_3 M_1^{-1} M_2 + M_4 \right) \Phi = \omega^2 M_5 \Phi.
%\end{equation}
%and is solved as a generalized eigenvalue problem. The choice of numerical fluxes which don't couple neighbour cells for $u$ leads to a sparse block diagonal matrix $M_1$ which can be easily inverted.\\
%For consistency using the solution $(u_0,\phi_0)$ of \eqref{eq:symm_mixed_form} we need to fulfill
%\begin{equation}
%\label{eq:mixed_form_consistency_conditions}
%\begin{cases}
%- \int_{\partial K} \widehat{u_0} \psi b \cdot n \dS= - \int_{\partial K} u_0 \psi b \cdot n \dS\qquad  & \forall \psi \in F_\Phi \\
%\int_{\partial K} \left( \widehat{\phi_0} -\phi_0^- \right)  v b \cdot n \dS =0 \qquad & \forall v \in F_U
%\end{cases}
%\end{equation}
To guarantee real eigenvalues, system matrices have to be symmetric which is ensured whenever the associated bilinear form is symmetric. The numerical fluxes are a slightly modified version of the local discontinuous Galerkin (LDG) fluxes proposed in \cite{outline_DG_cockburn_diffusion} and summarized in \cite[Table 3.1]{arnold} and given by
\begin{equation}
\label{eq:AD_LDG_variational_form_fluxes}
\widehat{u} = \fluxAvg{u}-\frac{\eta_S}{h_F} b\cdot \fluxJump{\phi}, \qquad \widehat{\phi} =  \fluxAvg{\phi}
\end{equation}
with the average $\fluxAvg{.}$ and jump $\fluxJump{.}$ on an interface $F$ defined as
\begin{equation}
\label{eq:AD_definition_flux_jump_averages}
\fluxAvg{f}  = \frac{1}{2}\left( f^K + f^{N_F(K)} \right)\quad ,\quad \fluxJump{f} = f^K \vn^K + f^{N_F(K)} \vn^{N_F(K)}
\end{equation}
where $\vn^K$ is the unit outer normal of cell $K$ and $\vn^{N_F(K)}$ is the unit outer normal of its neighbour  $N_F(K)$ sharing interface $F$. It holds $\vn^K = -\vn^{N_F(K)}$. We denote $h_F$ as the length of the respective cell edge at the interface. 
Further, $f^K$, $f^{N_F(K)}$ are the evaluations on $F$ from $K$ and from its neighbour respectively. 
%
%Further, $f^K$, $f^{N_F(K)}$ are defined as the limit on $F$ from $K$ and its neighbour respectively. For all $x \in F$ it holds
%\begin{align}
%f^K(x) & \coloneqq \lim\limits_{\substack{x_k \to x \\x_k \in K}} f(x_k)& , && f^{N_F(K)}(x) \coloneqq \lim\limits_{\substack{x_k \to x\\ x_k \in N_F(K)}} f(x_k) \ .
%\end{align}
We note that average and jump are well-defined as they are identical when viewed from each of the neighbouring cells of the interface. We further note that there are many possible choices for the fluxes. Setting $\eta_S = 0$ would lead to the Bassi-Rebay $1$ scheme \cite{outline_DG_bassi_rebay}, and decomposing $u$ into local and lifted gradients as done in the Bassi-Rebay $2$ scheme \cite{bassi_rebay_2} would result in a system where only direct neighbours are coupled.\\
%For the discretization using the locally aligned mesh, we choose this parameter to be $h_F = h_y$ for vertical interfaces and $h_F = h_x \sqrt{ 1+(b_2/b_1)^2}$ for $\vB$-aligned interfaces. The bilinear form then fulfills the required properties. 
Here, we choose the LDG fluxes \eqref{eq:AD_LDG_variational_form_fluxes} and insertion into \eqref{eq:variational_mixed_form} yields  $\forall K \in \mathcal{K}$
\begin{equation}
\label{eq:AD_LDG_variational_form_full}
\begin{dcases}
-\int_K u^K v^K \dV  + \int_K \vB \cdot \nabla \phi^K v^K \dV  &\\
\quad - \sum_{F \in \partial K} \left( \int_{F} \frac{1}{2} \vB \cdot \left( \phi^K \vn^K + \phi^{N_F(K)} \vn^{N_F(K)} \right)  v^K  \dS \right)  =0 &  \forall v_K \in V_{K,U}\\
\ \\
- \sum_{F \in \partial K}  \left( \int_F \frac{1}{2}\left( u^K + u^{N_F(K)}  \right) \psi^K \vB \cdot \vn^K \dS \right) &\\
\quad + \sum_{F \in \partial K} \left(  \int_{F}  \frac{\eta_S}{h_F} \vB \cdot \left( \phi^K \vn^K + \phi^{N_F(K)} \vn^{N_F(K)} \right)  \psi^K \vB \cdot \vn^K \dS \right) &\\
\quad + \int_K u^K \vB \cdot \nabla \psi^K \dV  = \omega^2 \int_K \alpha \phi^K \psi^K \dV  & \forall \psi^K \in V_{K,\Phi}
\end{dcases}
\end{equation}
To retrieve the global variational formulation, we sum \eqref{eq:AD_LDG_variational_form_full} over all over all elements~$K$. For this we establish the following Lemmata by defining $\mathcal{F}$ as the space of all interfaces.
\begin{Lemma}
	\label{thm:AD_simplification_summation_flux_arbitrary}
	Let $\widehat{\phi}$ be an arbitrary flux, $\psi^K \in V_{K,\Phi}$ locally defined in $K$ and $\vv \in \C^2$ a not necessarily constant vector. Then
	\begin{equation}
	\sum_{K \in \mathcal{K}} \sum_{F \in \partial K} \left( \int_{F} \widehat{\phi} \psi^K \vv \cdot \vn^K \dS\right) = \sum_{F \in \mathcal{F}} \left( \int_{F} \widehat{\phi} \vv \cdot \fluxJump{\psi} \dS\right) \ . 
	\end{equation}
\end{Lemma}
\begin{Lemma}
	\label{thm:AD_simplification_summation_flux_avg}
	Let $\widehat{\phi} = \fluxAvg{\phi}$, $\psi^K \in V_{K,\Phi}$ locally defined in $K$ and $\vv \in \C^2$ a not necessarily constant vector. Then
	\begin{equation}
	- \sum_{K \in \mathcal{K}} \sum_{F \in \partial K} \left( \int_{F} \left( \widehat{\phi} - \phi^K \right) \psi^K \vv \cdot \vn^K \dS\right) = \sum_{F \in \mathcal{F}} \left( \int_{F} \vv \cdot \fluxJump{\phi} \fluxAvg{\psi} \dS\right) \ . 
	\end{equation}
\end{Lemma}
%\Proof{See Appendix~\ref{sec:Appendix_proof_lemma_2}.}
\Proof{ \textblue{For both Lemma~\ref{thm:AD_simplification_summation_flux_arbitrary} and \ref{thm:AD_simplification_summation_flux_avg}: When summing over all cells, each interface $F$ is considered twice in total, once for a cell $K$ and once for its unique neighbour $N_F(K)$. The result is then obtained by straightforward calculation. }
}
Using Lemma~\ref{thm:AD_simplification_summation_flux_avg} on the surface term of the first equation of \eqref{eq:AD_LDG_variational_form_full} and summing over all cells $K \in \mathcal{K}$, we obtain
\begin{equation}
\label{eq:AD_LDG_summed_integrals_start}
\sum_{K \in \mathcal{K}}  \left( - \int_K u^K v^K \dV  + \int_K \vB \cdot \nabla \phi^K v^K \dV \right) - \sum_{F \in \mathcal{F}}  \int_{F} \vB \cdot \fluxJump{\phi} \fluxAvg{v} \dS =0
\end{equation}
whereas the second equation of \eqref{eq:AD_LDG_variational_form_full} using Lemma~\ref{thm:AD_simplification_summation_flux_arbitrary} yields
\begin{equation}
\label{eq:AD_LDG_summed_integrals_end}
\begin{aligned}
- &\sum_{F \in \mathcal{F}} \left( \int_{F} \fluxAvg{u} \vB \cdot \fluxJump{\psi} \dS \right) + \sum_{F \in \mathcal{F}} \left( \int_{F}  \frac{\eta_S}{h_F} \vB \cdot \fluxJump{\phi} \vB \cdot \fluxJump{\psi}  \dS \right) \\
& + \sum_{K \in \mathcal{K}}  \left( \int_K u^K \vB \cdot \nabla \psi^K \dV \right) = \omega^2 \sum_{K \in \mathcal{K}}   \int_K \alpha \phi^K \psi^K \dV  \ .
\end{aligned}
\end{equation}
The method is consistent as for a continuous solution $\left( \phi_0, u_0,  \omega_0^2 \right)$ it holds $\fluxAvg{\phi_0} = \phi_0, \fluxAvg{u_0} = u_0$ and $\fluxJump{\phi_0} = 0$. Insertion in \eqref{eq:AD_LDG_summed_integrals_start} and \eqref{eq:AD_LDG_summed_integrals_end} and one integration by parts in \eqref{eq:AD_LDG_summed_integrals_end} allows us to retrieve \eqref{eq:AD_mixed_form_var_form_general}.\\
For building the system matrices and to prove its symmetry, we now associate the matrix components as follows\\
%\begin{align}
%\label{eq:matrix_integral_associations_auv}
%M_{U V} &\leftrightarrow \sum_K \int_K u^K v^K \dV  \\
%\label{eq:matrix_integral_associations_aphiv}
%A_{\Phi V} & \leftrightarrow \sum_K \int_K \vB \cdot \nabla \phi^K v^K \dV \\ 
%\label{eq:matrix_integral_associations_bphiv}
%B_{\Phi V}&\leftrightarrow \sum_{F}  \int_{F} \frac{1}{2}\left( \phi^{K_F^-} \vB \cdot \vn^{K_F^-} +  \phi^{K_F^+}  \vB \cdot \vn^{K_F^+} \right)  \left( v^{K_F^-} +v^{K_F^+} \right)  \dS \\
%\label{eq:matrix_integral_associations_bupsi}
%B_{U \Psi} &\leftrightarrow \sum_F \int_{F} \frac{1}{2}\left( u^{K_F^-} + u^{K_F^+} \right) \left( \psi^{K_F^-} \vB \cdot \vn^{K_F^-} + \psi^{K_F^+} \vB \cdot \vn^{K_F^+} \right)  \dS \\
%\label{eq:matrix_integral_associations_bphipsi}
%B_{\Phi \Psi} &\leftrightarrow \sum_F  \int_{F}  \frac{\eta_S}{h_F} \left( \phi^{K_F^-} \vB \cdot \vn^{K_F^-} + \phi^{K_F^+} \vB\cdot \vn^{K_F^+}\right) \left(  \psi^{K_F^-} \vB \cdot \vn^{K_F^-}  + \psi^{K_F^+} \vB \cdot \vn^{K_F^+} \right) \dS \\
%\label{eq:matrix_integral_associations_aupsi}
%A_{U \Psi} & \leftrightarrow \sum_K  \int_K u^K \vB \cdot \nabla \psi^K \dV   \\
%\label{eq:matrix_integral_associations_aphipsi}
%M_{\Phi \Psi} &\leftrightarrow \sum_K  \int_K \phi^K \psi^K \dV   
%\end{align}\\
\\
\begin{minipage}{0.48\textwidth}
	\begin{align}
	\label{eq:AD_LDG_system_matrices_associations_auv}
	M_{U V} &\leftrightarrow  \sum_{K \in \mathcal{K}}   \int_K u^K v^K \dV   \\
	\label{eq:AD_LDG_system_matrices_associations_aphiv}
	A_{\Phi V} & \leftrightarrow \sum_{K \in \mathcal{K}}  \int_K \vB \cdot \nabla \phi^K v^K \dV   \\ 
	\label{eq:AD_LDG_system_matrices_associations_bphiv}
	B_{\Phi V}&\leftrightarrow \sum_{F \in \mathcal{F}}  \int_{F} \vB \cdot \fluxJump{\phi} \fluxAvg{v}   \dS  \\
	\label{eq:AD_LDG_system_matrices_associations_bupsi}
	B_{U \Psi} &\leftrightarrow \sum_{F \in \mathcal{F}}  \int_{F} \fluxAvg{u} \vB \cdot \fluxJump{\psi} \dS   
	\end{align}
\end{minipage}
\begin{minipage}{0.48\textwidth}
	\begin{align}
	\label{eq:AD_LDG_system_matrices_associations_bphipsi}
	B_{\Phi \Psi} &\leftrightarrow \sum_{F \in \mathcal{F}}  \int_{F}  \frac{\eta_S}{h_F} \vB \cdot \fluxJump{\phi} \vB \cdot \fluxJump{\psi}  \dS  \\
	\label{eq:AD_LDG_system_matrices_associations_aupsi}
	A_{U \Psi} & \leftrightarrow \sum_{K \in \mathcal{K}}   \int_K u^K \vB \cdot \nabla \psi^K \dV   \\
	\label{eq:AD_LDG_system_matricesl_associations_aphipsi}
	M_{\Phi \Psi} &\leftrightarrow \sum_{K \in \mathcal{K}}   \int_K \alpha \phi^K \psi^K \dV    
	\end{align}
\end{minipage}\\
\ \\
The system then writes
\begin{equation}
\label{eq:mixed_form_system_full}
\begin{pmatrix}
-M_{U V} &   A_{\Phi V} - B_{\Phi V}  \\
A_{U \Psi} - B_{U \Psi} & B_{\Phi \Psi}
\end{pmatrix}
\begin{pmatrix}
U\\
\Phi
\end{pmatrix}
= \omega^2
\begin{pmatrix}
0 & 0 \\
0 & M_{\Phi \Psi}
\end{pmatrix}
\begin{pmatrix}
U\\
\Phi
\end{pmatrix}.
\end{equation}
Choosing the same basis for test and trial functions, it is obvious that the mass matrices $M_{U V}$ and $M_{\Phi \Psi}$ and the penalization matrix $B_{\Phi \Psi}$ are symmetric. Furthermore, it holds $A_{\Phi V} = A_{U \Psi}^\top$ and $B_{\Phi V} = B_{U \Psi}^\top$ as $\phi, \psi$ and $u, v$ are interchangeable between \eqref{eq:AD_LDG_system_matrices_associations_aphiv} and \eqref{eq:AD_LDG_system_matrices_associations_aupsi} as well as \eqref{eq:AD_LDG_system_matrices_associations_bphiv} and \eqref{eq:AD_LDG_system_matrices_associations_bupsi}.\\
\eqref{eq:mixed_form_system_full} is reduced to the generalized eigenvalue problem
\begin{equation}
\begin{aligned}
\label{eq:reduced_system_matrix}
	 A &\coloneqq \left( \left(A_{U \Psi} - B_{U \Psi}\right) M_{U V}^{-1}  \left(A_{U \Psi}^\top - B_{U \Psi}^\top\right) + B_{\Phi \Psi} \right) \\
	 A\Phi &  = \omega^2 M_{\Phi \Psi} \Phi.
\end{aligned}
\end{equation} 
We note that the matrix $A$ is symmetric and the mass matrix $M_{UV}$ only has element-local contributions and hence is an easily invertible block diagonal matrix. Also note that \eqref{eq:reduced_system_matrix} can be written as a standard symmetric eigenvalue problem using the root of the mass matrix $M_{\Phi \Psi}$ 
\begin{equation}
	 M_{\Phi \Psi}^{-\frac{1}{2}} A M_{\Phi \Psi}^{-\frac{1}{2}} \left( M_{\Phi \Psi}^\frac{1}{2} \Phi \right)   = \omega^2 \left(  M_{\Phi \Psi}^\frac{1}{2} \Phi \right) .
\end{equation}
The eigenvalue problem is solved with the \texttt{FEAST} library \cite{feast_documentation} which allows to solve for a specific region of the spectrum. We also use the \texttt{MUMPS} library \cite{MUMPS_1,MUMPS_2} for performing matrix multiplications and factorizations as well as the solution of linear systems needed by the \texttt{FEAST} eigenvalue solver.

\section{Numerical results for constant coefficients}
\label{sec:results}

As analytic solutions are available for $\alpha, \beta \equiv 1$, we compare the results of the proposed locally field-aligned discontinuous Galerkin method (ADG) with the exact eigenvalues deduced in Section~\ref{sec:eigenfunctions_gradients}. Therefore, the discrete eigenfunctions need to be associated to the exact eigenfunctions in a postprocessing step which is described in Section~\ref{sec:eigenvector_postprocessing}. For the assessment of ADG, we first define a reference case for one specific choice of constant $\vB$ in Section~\ref{sec:reference_case}. We then examine the impact of the local alignment of mesh and basis in Section~\ref{sec:impact_of_alignment}. In Section~\ref{sec:dof_ratio}, we investigate different choices for the distribution of parallel and perpendicular resolution. We investigate the convergence rate of ADG in  Section~\ref{sec:convergence}.\\
The total number of degrees of freedom is given by
\begin{equation}
\dof = \dof_\parallel \dof_\perp = (p_\xi+1) N_x (p_\eta+1)  N_y.
\end{equation}
Throughout this chapter, we present relative errors whenever we speak of errors. If the exact eigenvalue of a mode is zero, we use absolute errors instead. This is particularly the case for the constant mode.\\
%For comparisons, the quality measure for the method deals with the size of the biggest error of the eigenvalues associated to a restricted space of Fourier modes.\\
Usually the biggest error occurs for eigenvalues which are associated to eigenmodes of the highest considered mode number. As we will see, this is not always the case for ADG. Due to the separation of parallel and perpendicular resolution and the analysis of Section~\ref{sec:eigenfunctions_gradients}, we expect the error of eigenvalues to scale with the size of the eigenvalue and the mode number of the associated eigenfunction.\\

\subsection{Eigenvector post-processing}
\label{sec:eigenvector_postprocessing}

For each eigenvalue $\omega_j$, a discrete eigenvector $\left\{ \Phi_k^j \right\}_{k=1}^{\dof}$ is found. Each degree of freedom is associated to a discontinuous Galerkin basis function $\phi^k$.\\
We compute the projection of each basis function onto the set of Fourier modes $ \abs{m} \leq m_\text{max}$, $ \abs{n}  \leq n_\text{max}$
\begin{equation}
\tilde{\phi}_{m,n}^k = \int_\Omega \phi^k \phi_{m,n} \dV \ . 
\end{equation}
We then associate the eigenvalue $\omega_j$ to the mode number with the maximal amplitude given by
\begin{equation}
\underset{\abs{m}\leq m_\text{max},\abs{n} \leq n_\text{max}}{\argmax} \abs {\sum_{k = 1}^{\dof} \vPhi_k^j  \tilde{\phi}_{m,n}^k} \ . 
\end{equation}

%We associate each discrete eigenvector $\vPhi^j$ to a mode number by projection onto the space of exact eigenmodes $\phi_{m,n}$ \eqref{eq:exact_eigenfunctions} with limited mode numbers $ \left| m \right| \leq m_\text{max}$, $ \left| n \right|  \leq n_\text{max}$ yielding a spectrum
%\begin{equation}
%\vsigma_{m,n}^j = \int_\Omega \vPhi^j \phi_{m,n} \dV.
%\end{equation}
%The mode number with the maximal amplitude given by 
%\begin{equation}
%\underset{m,n}{\argmax} \left| \sigma_{m,n}^j \right|
%\end{equation}
%is associated to the discrete eigenvector $\vPhi_j$.

\subsection{Reference case}
\label{sec:reference_case}

\begin{figure}[!htb]
	\centering
	\includegraphics[width=0.6\textwidth]{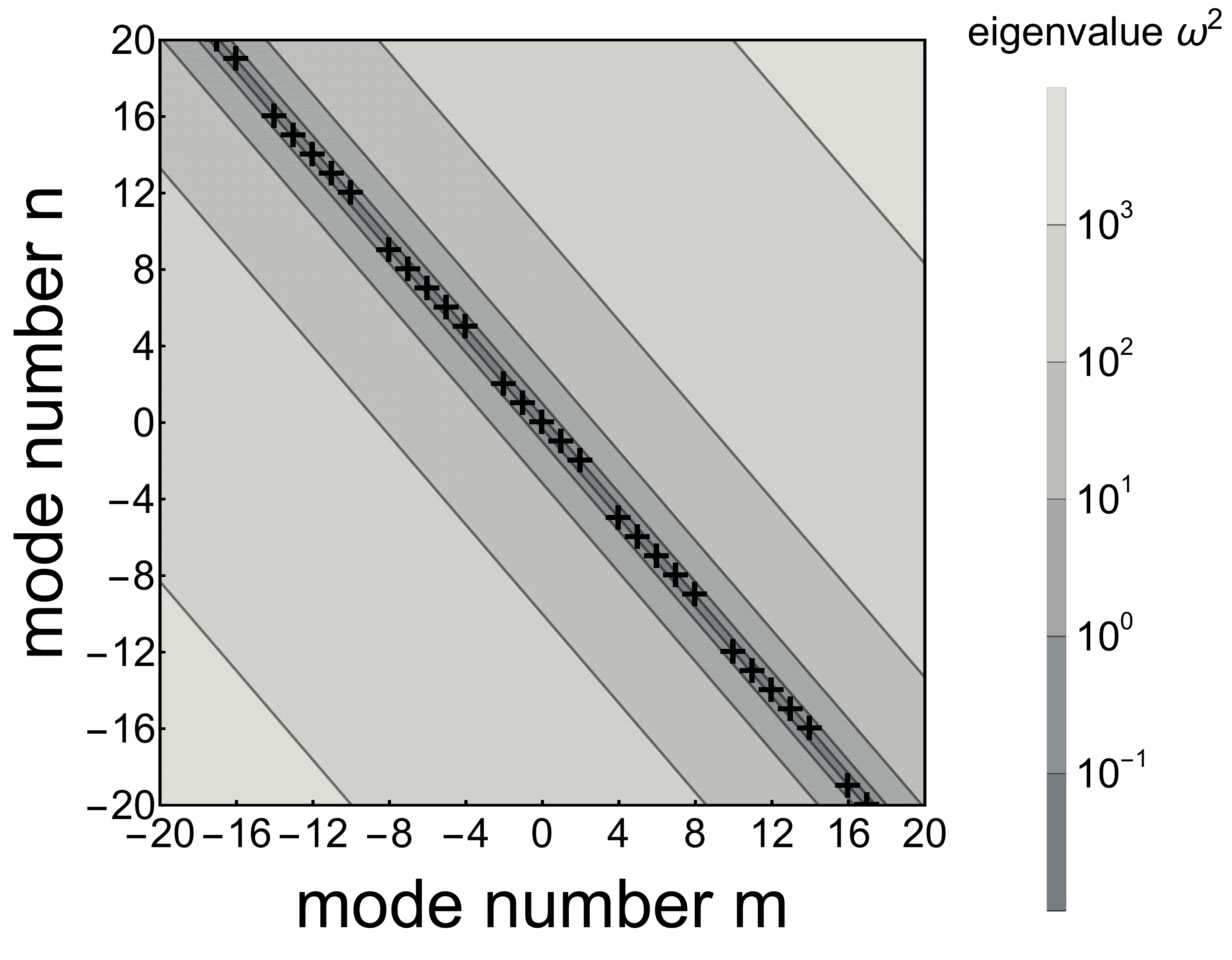}
	\caption{Contour plot of the size of the exact eigenvalues for the reference case with associated modes with maximal frequency $m_{\text{max}} = n_{\text{max}} = 20$. \textred{The black crosses indicate the discrete mode numbers $(m,n)$ that lie within the band of associated eigenvalues $\omega^2<0.2$.}} 
	\label{fig:modes_exact_eigenvalues}
\end{figure}

For the purpose of evaluating the properties of ADG, we define a reference test case by \textred{$\vb = \left( \iota, 1 \right)^\top = (  1.165939761, 1)^\top $ } . This value for $\vb$ is chosen to avoid a low rational number for $\iota$ such that all considered eigenvalues differ from zero except for the constant mode. Further, it yields a locally aligned mesh with low shear and non-conforming interfaces of different lengths.\\
In the reference case, we consider mode numbers up to $m_\text{max} = n_\text{max} = 20$. Figure~\ref{fig:modes_exact_eigenvalues} shows the distribution of eigenvalues on a logarithmic scale for the given selection of modes. We observe that small eigenvalues gather perpendicularly to $\vb$.
%\begin{figure}[!htb]
%	\centering
%	\includegraphics[width=0.5\textwidth]{Figures/modes_exact_eigenvalues.eps}
%	\caption{Contour plot of the size of the exact eigenvalues for the reference case with associated modes with maximal frequency $m_{\text{max}} = n_{\text{max}} = 20$. The white dashed line indicates the direction perpendicular to $\vB$.} 
%	\label{fig:modes_exact_eigenvalues}
%\end{figure}\\
\\
As modes with small parallel gradient are of interest, we aim to resolve all modes with associated eigenvalue $\omega_{m,n}^2 \leq 0.2$. The stabilization parameter $\eta_S$ in \eqref{eq:AD_LDG_variational_form_fluxes} is chosen as $6$.\\
%The exact eigenvalues and associated eigenmodes are summarized in Table~\ref{tab:exact_eigenvalues}.\\
%\begin{center}
%	\begin{tabular}{cc}
%		eigenvalue & eigenmode  \\
%		\hline
%		$0$ & $(0,0)$ \\
%		$0.0275360047591036$ & $(1,-1)$ \\ 
%		$0.1101440190364146$ & $(2,-2)$ \\ 
%		$0.1130579766145697$ & $(4,-5)$ \\  
%		$0.0290024945637303$ & $(5,-6)$ \\
%		$0.0000190220310983$ & $(6,-7)$ \\
%		$0.0261075590166736$ & $(7,-8)$ \\
%		$0.1072681055204561$ & $(8,-9)$
%	\end{tabular}
%	\label{tab:exact_eigenvalues}
%	\captionof{table}{Exact eigenvalues and eigenmodes for the reference case.}	
%\end{center}
%We remark that system matrices are stored as lower triangular sparse matrices. Therefore, we investigate ADG regarding the numbers of non-zeroes of the lower triangle of the left hand side system matrices in relation to the total number of matrix entries and abbreviate this percentage by nnzA.
We remark that system matrices are stored as lower triangular sparse matrices. nnzA is the percentage of non-zeroes entries of the system matrix $A$.

%Therefore, we investigate ADG regarding the numbers of non-zeroes of the lower triangle of the left hand side system matrices in relation to the total number of matrix entries and abbreviate this percentage by nnzA.

\subsection{Impact of the local alignment}
\label{sec:impact_of_alignment}

\begin{figure}[!htb]
	\subfigure[cartesian mesh] {\includegraphics[width=0.48\textwidth]{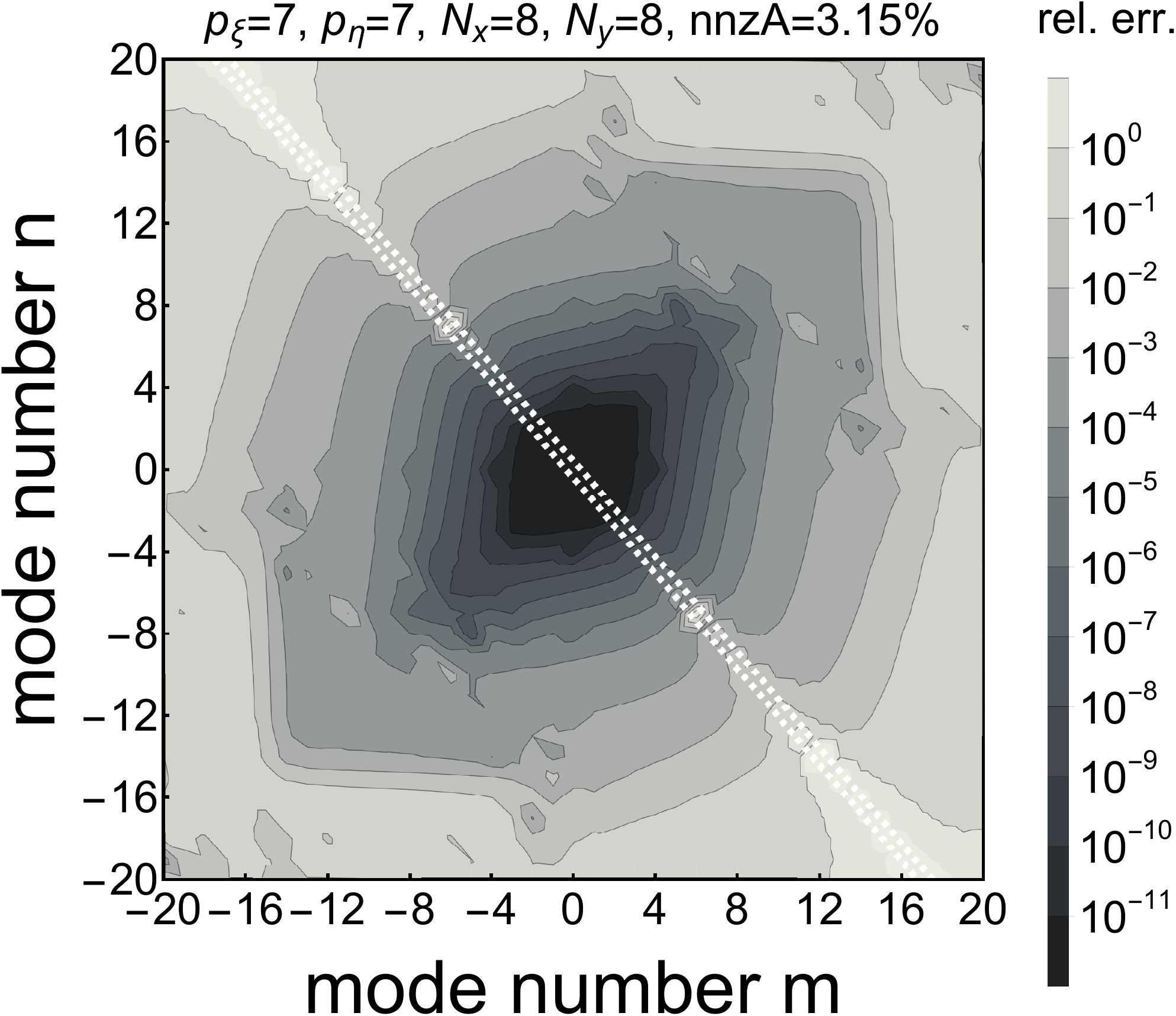}}
	\subfigure[locally aligned mesh] {\includegraphics[width=0.48\textwidth]{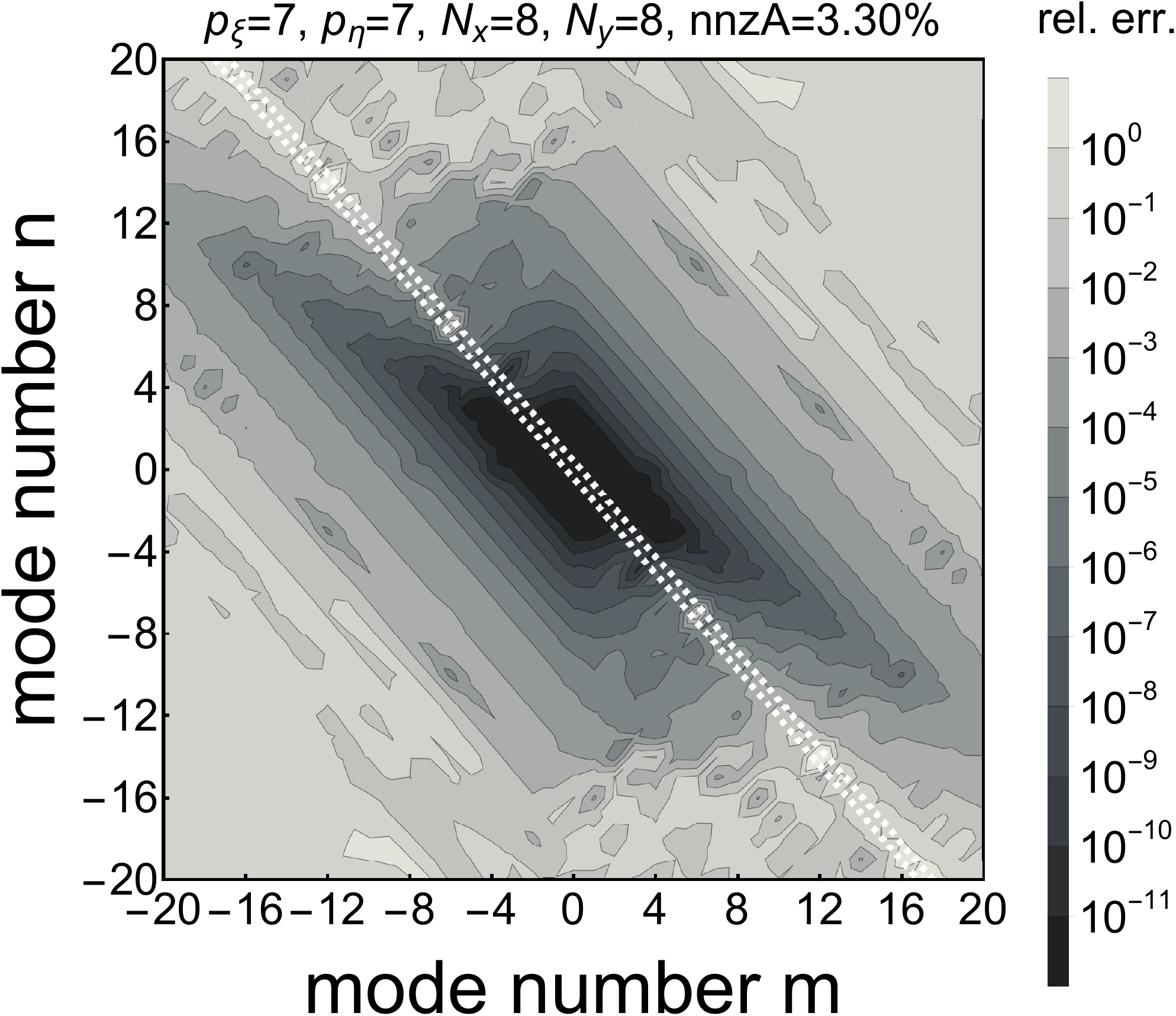}}
	\caption{Errors of a non-aligned cartesian mesh DG method and the locally aligned mesh of ADG, in the reference case with $\dof = 2^{12}$. In between the white dashed lines resides the band of modes with eigenvalues $\omega^2 \leq \omega_{\operatorname{max}}^2 = 0.2$.} 
	\label{fig:impact_alignment_all_modes}
\end{figure}

\begin{figure}[!htb]
	\centering
	\includegraphics[width=0.9\textwidth]{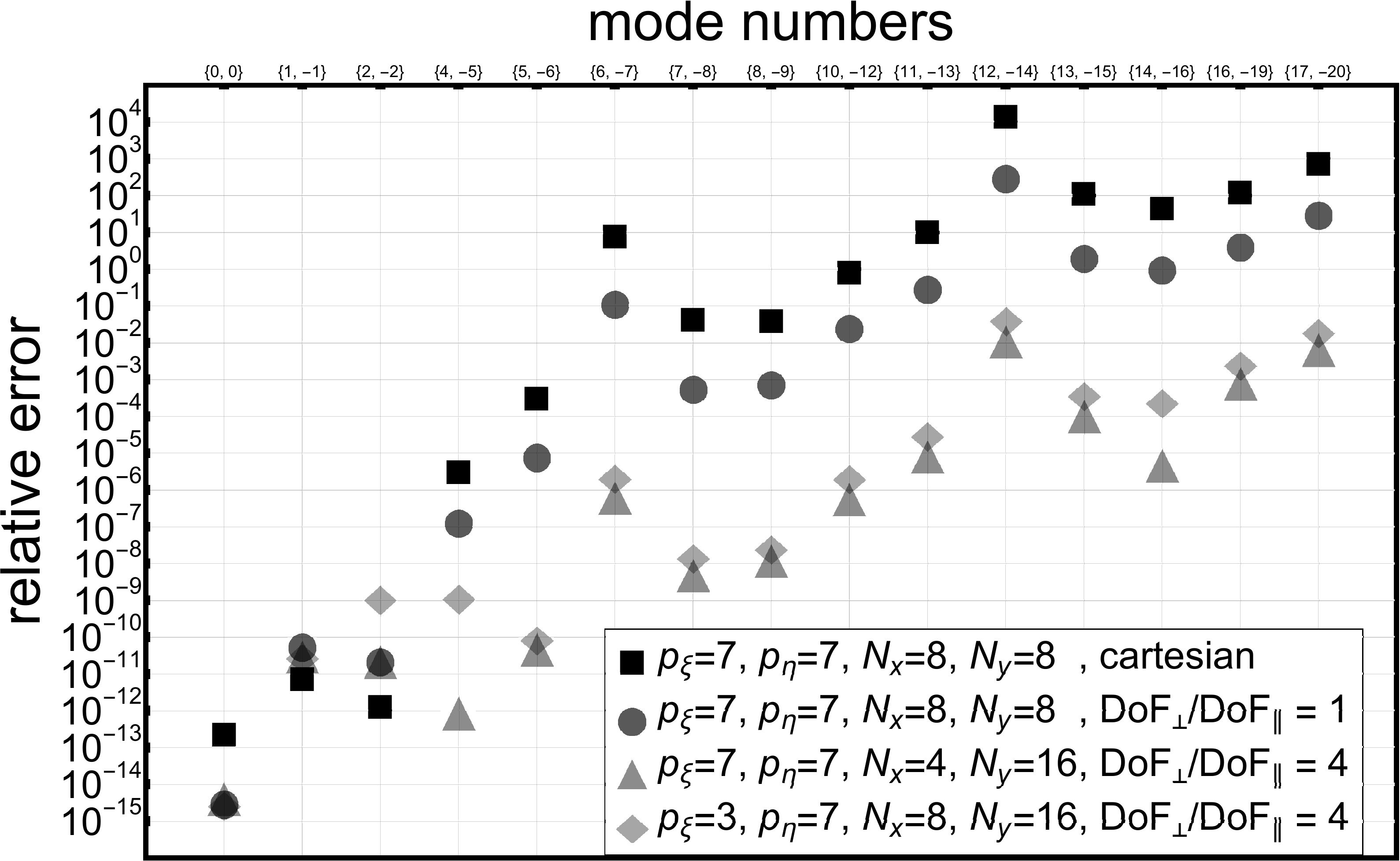}
	\caption{Comparison of errors of a non-aligned cartesian case with ADG in the reference case with $\dof = 2^{12}$ for different distributions of resolution.} 
	\label{fig:impact_alignment_comparison}
\end{figure}

In this section, we compare a non-aligned DG method operating on a cartesian mesh with the locally field-aligned mesh of ADG for the same number of degrees of freedom. We choose the polynomial degrees $p_\xi = 7 = p_\eta$ and the mesh resolution $N_x = 8 = N_y$.\\
In Figure~\ref{fig:impact_alignment_all_modes}, we first show the eigenvalue errors for all Fourier modes $\abs{m}, \abs{n} \leq 20$ as contour line plots. For the cartesian mesh used in Figure~\ref{fig:impact_alignment_all_modes}(a), we observe that the error increases for higher mode numbers. Similar behaviour is found in Figure~\ref{fig:impact_alignment_all_modes}(b) as well, but additionally a correlation of the error to the magnitude of the eigenvalue is introduced. The well resolved region is tilted towards the direction perpendicular to $\vb$, which is the region where small eigenvalues reside as indicated by the white dashed line in Figure~\ref{fig:modes_exact_eigenvalues}.\\
In Figure~\ref{fig:impact_alignment_comparison}, we take a closer look at the band of modes with eigenvalues $\omega^2 \leq \omega_\text{max}^2$. We observe that the errors of eigenvalues with large mode numbers modes are smaller by $1.5$ to $2$ orders of magnitude for ADG (round markers) in comparison to a non-aligned discontinuous Galerkin method (square makers). Thus, we conclude that aligning the mesh yields a significant accuracy improvement within the selected mode band.
\\
\\

\subsection{Distribution of resolution}
\label{sec:dof_ratio}

When comparing the band of modes for $\omega^2 \leq \omega_\text{max}^2$ in Figure~\ref{fig:impact_alignment_all_modes} and the distribution of exact eigenvalues in Figure~\ref{fig:modes_exact_eigenvalues}, we observe that many eigenvalues of mode numbers outside this band are overresolved. \textblue{As shown in Section~\ref{sec:eigenfunctions_gradients}, the  eigenfunctions with small eigenvalues have a small parallel gradient.} Thus, we aim to distribute the resolution $\dof_\parallel$ and $\dof_\perp$ of the method such that $\dof_\perp > \dof_\parallel$.
\begin{figure}[!htb]
	\subfigure[$p_\xi = 7$] {\includegraphics[width=0.47\textwidth]{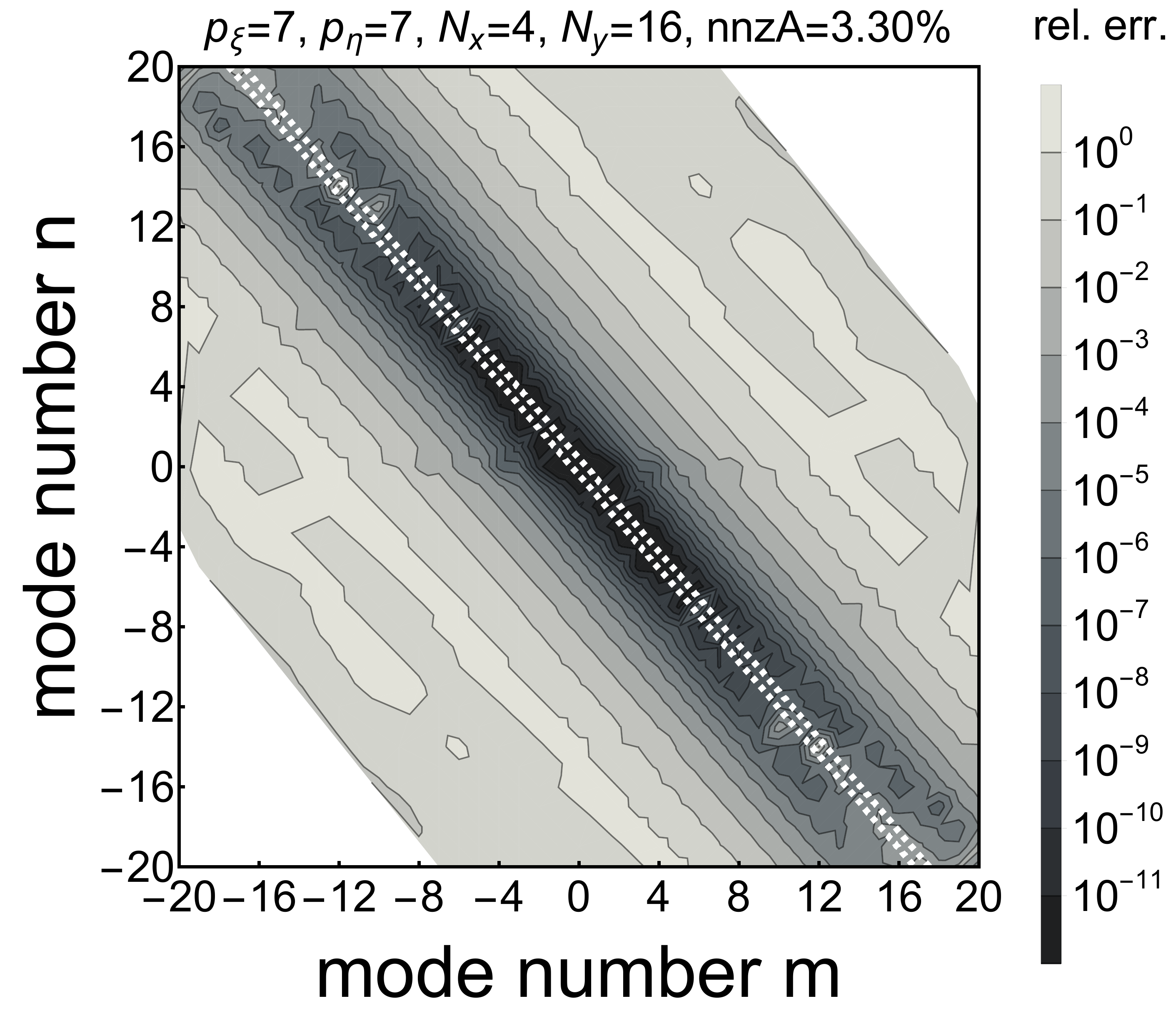}}
	\subfigure[$p_\xi = 3$] {\includegraphics[width=0.47\textwidth]{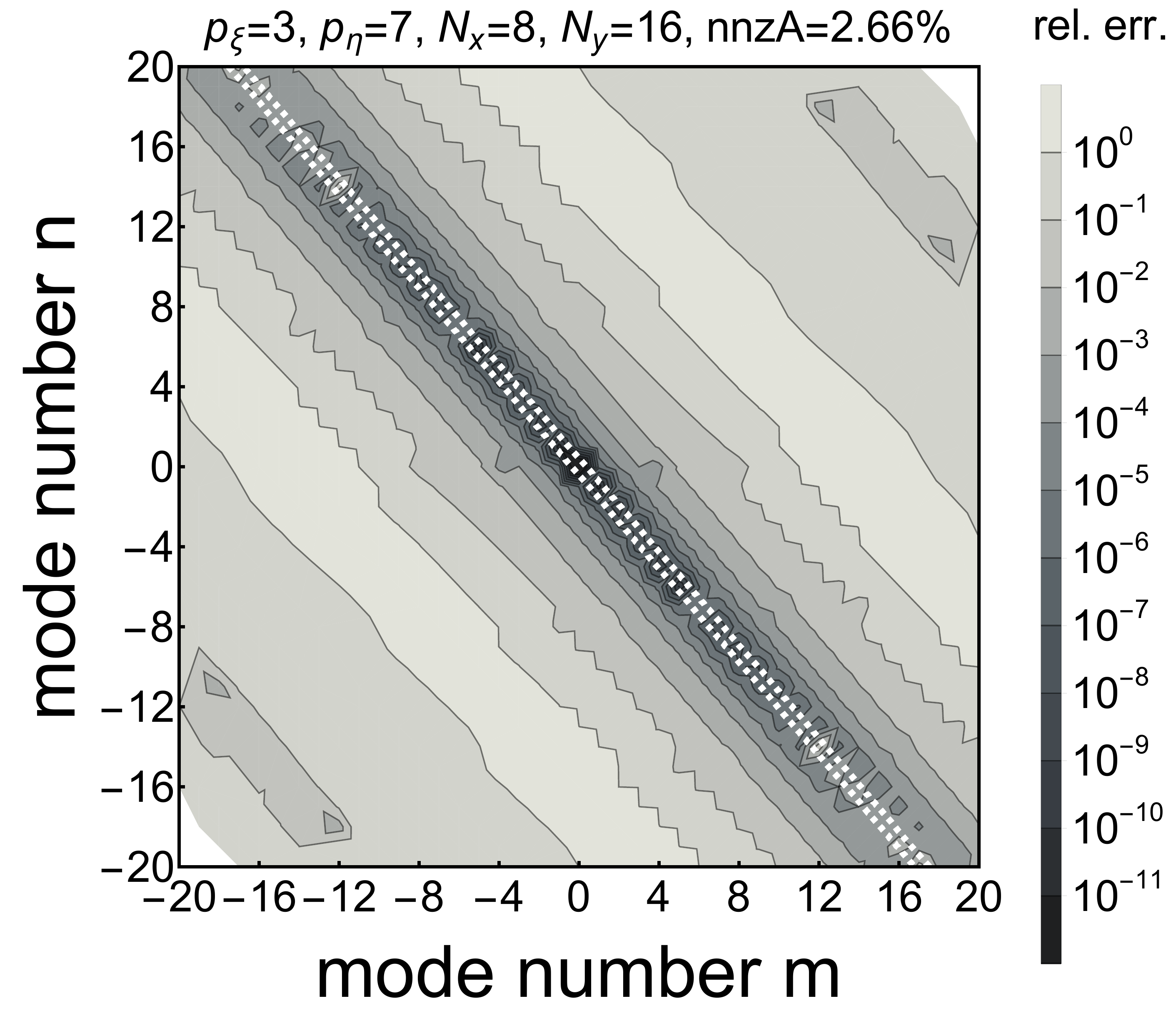}}
	\caption{Errors of ADG for high and low degree $p_\xi$ but same parallel resolution $\dof_\parallel$ with $\dof = 2^{12}$, in the reference case. In between the white dashed lines resides the band of modes with eigenvalues $\omega^2 \leq \omega_{\operatorname{max}}^2 = 0.2$.} 
	\label{fig:distribution_resolution_all_modes}
\end{figure}
First, we keep the total number of cells constant and change the cell distribution by refining $N_y$ and coarsening $N_x$. The effects are shown in Figure~\ref{fig:distribution_resolution_all_modes}(a). Comparing these results to  to Figure~\ref{fig:impact_alignment_all_modes}(b), we observe that the region of well-resolved eigenvalues with errors smaller than $10^{-4}$ extends into regions of larger mode numbers and gathers narrower around the interesting band of modes.\\
The errors of the mode band are also plotted in Figure~\ref{fig:impact_alignment_comparison}. We observe that changing the cell distribution yields an increase in accuracy of $4$ to $5$ orders of magnitude for mode numbers larger than $4$ when comparing round markers with triangle markers. The comparison with a cartesian mesh (square markers) yields $5.5$ to $7$ orders of magnitude in total by aligning the mesh and distribute its resolution such that $\dof_\perp/\dof_\parallel = 4$ for the same total resolution $\dof = 2^{12}$.\\
We can also adapt the ratio $\dof_\perp/\dof_\parallel$ by modifying the degree of the basis functions. Figure~\ref{fig:distribution_resolution_all_modes}(b) shows a configuration with $p_\xi = 3$ and $\dof_\perp/\dof_\parallel=4$. We observe that the well resolved region is thinner whereas the errors related to the mode number are of the same magnitude in comparison to Figure~\ref{fig:distribution_resolution_all_modes}(a). The accuracy of results within the interesting band of modes is marginally affected as confirmed by Figure~\ref{fig:impact_alignment_comparison} when comparing $p_\xi =7$ and $p_\xi = 3$ (triangle markers and diamond-shaped markers) with $\dof_\perp/\dof_\parallel = 4$. However, a lower degree $p_\xi$ increases the sparsity of system matrices from $3.30\%$ to $2.66\%$, so we trade some accuracy for higher sparsity.

\subsection{Numerically observed convergence}
\label{sec:convergence}

For examining the convergence behaviour of ADG, we trace the maximal error inside the band of modes $\omega^2\leq \omega_{\operatorname{max}}^2$ with mode numbers up to $20$. We consider a parallel degree of $p_\xi = 3$ and  perpendicular degree of $p_\eta = 3,7$ and refine the mesh in $N_x, N_y$ simultaneously by doubling the number of cells in each direction. We consider this for configurations with $\dof_\perp / \dof_\parallel = 4,8$. 
\begin{figure}[!htb]
	\centering
	\includegraphics[width=0.9\textwidth]{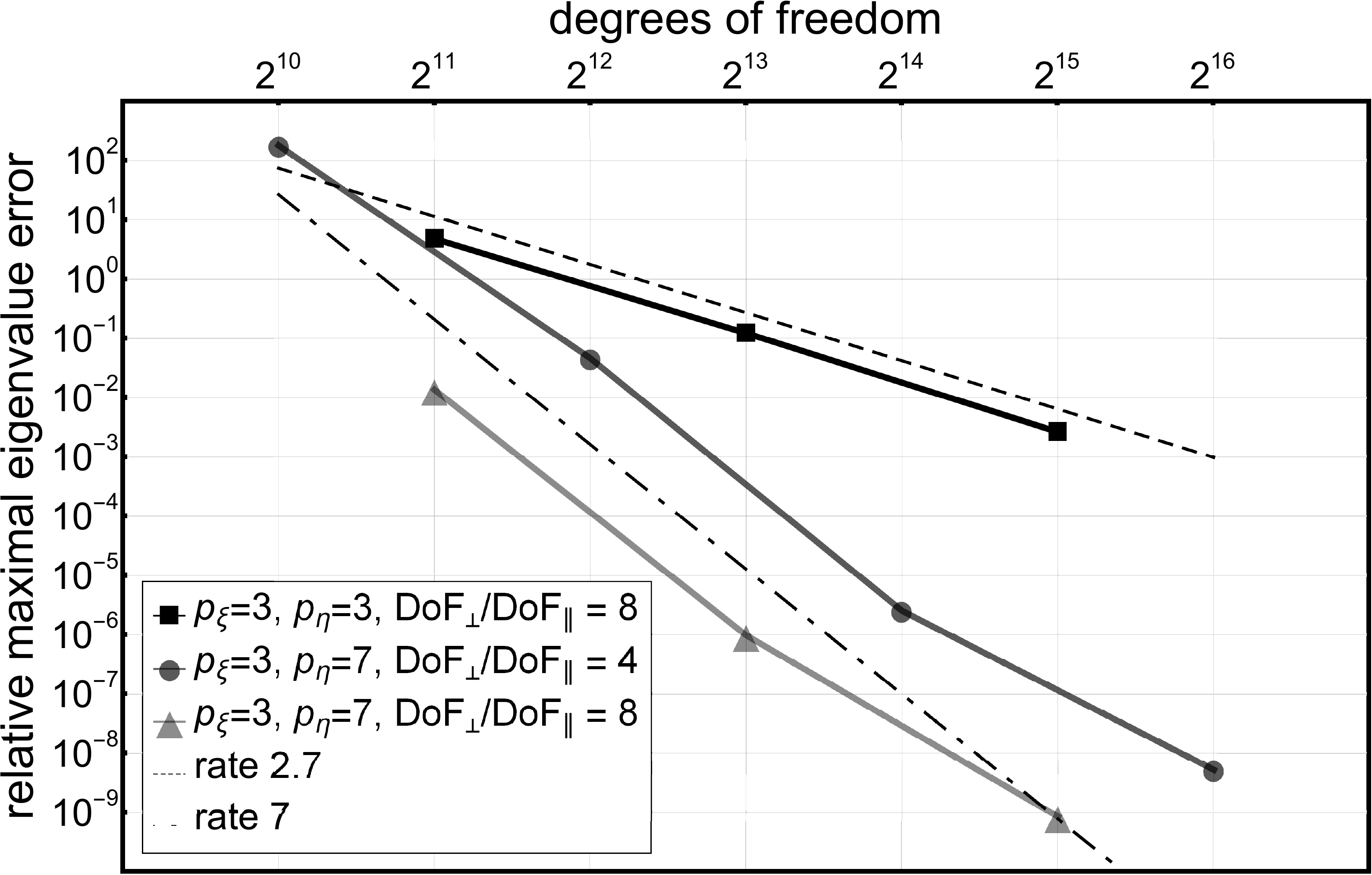}
	\caption{Convergence of the maximal error among all considered eigenmodes in the reference case. Different mesh configurations and refinement strategies are considered. The numerically observed convergence rate is plotted as dashed lines.}
	\label{fig:convergence_ADG}
\end{figure}\\
Figure~\ref{fig:convergence_ADG} shows the convergence of these configurations. For the same resolution, $\dof_\perp / \dof_\parallel = 8$ yields an improvement of $2.5$ to $3$ orders of magnitude compared to $\dof_\perp / \dof_\parallel = 4$. The numerically observed rate of convergence of the maximal error in the reference configuration is of the order of the perpendicular degree $p_\eta$ and approximately given by
\begin{equation}
\mathcal{O} \left( \dof^{-p_\eta}\right) \ .
\end{equation}
%We further consider a refinement strategy using fixed $\dof_\parallel$ which exclusively refines the mesh in $N_y$ by doubling the number of cells. Figure~\ref{fig:convergence_ADG} shows, that we obtain rapid convergence in this case (triangle markers). This process lasts until the given parallel resolution is too coarse to obtain better results. The numerically observed rate of convergence is doubled. 
As we consider the convergence of multiple eigenvalues at once, we leave the convergence rates as a bare observation rather than stating it as a general property of ADG.

\section{Numerical results for a three-dimensional MHD equilibrium}
\label{sec:ADM_MHD}

\begin{figure}[!htb]
	\centering
	\subfigure[$ \alpha\left(\vx\right)\quad\quad$]{\includegraphics[width=0.45\textwidth]{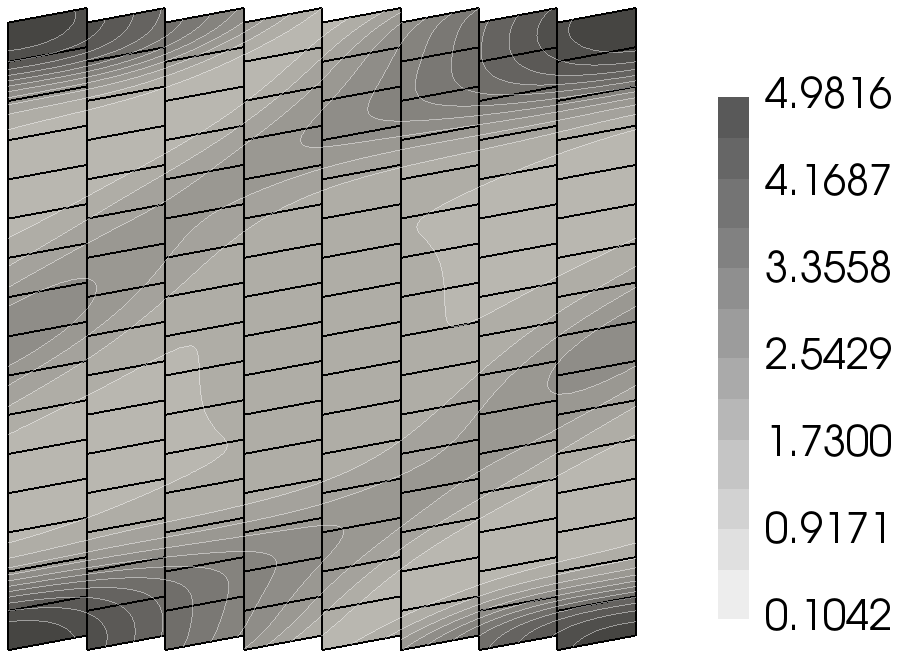}
}
\subfigure[$ \beta\left(\vx\right)\quad\quad$]{\includegraphics[width=0.45\textwidth]{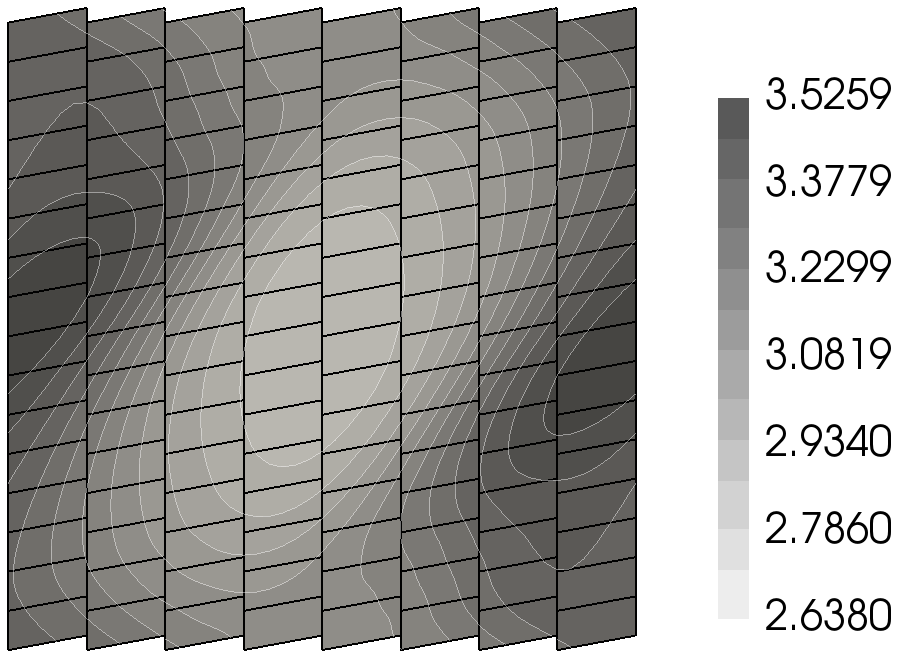}
} 
	\caption{\textred{Distribution of the variable coefficients for the flux surface $s=0.8$ of the W7-X case.}}
	\label{fig:varcoefs}
\end{figure}

\textblue{In this section, we investigate the ADG method for the anisotropic wave equation with \emph{variable coefficients}. The two-dimensional periodic domain is mapped to a flux surface of a three-dimensional  MHD equilibrium, with two periodic angles $\theta,\varphi$, chosen such that the magnetic field direction is constant in the logical domain. 
The coefficients $\alpha \left( \vx \right)$ and $\beta \left( \vx \right)$ are then computed from the mapping and the magnetic field of that flux surface, see \cite[Section~2.6]{mydiss} for details. As a test case, we use the \texttt{VMEC}-equilibrium of the W7-X high-mirror case \cite[Table IV]{vmec_equilibrium}.\\
We solve the two-dimensional eigenvalue problem on a sequence of nested flux surfaces, which are parameterized by a normalized flux surface coordinate $s \in [0,1]$ with $0$ being the magnetic axis and $1$ being the outermost flux surface.} \textred{ In the considered W7-X case, the magnetic field is given by $\vb = \left( \iota (s) , 1\right)^\top$ with $\iota(s)=0.85931(1-s)+0.93972s$}. 
%The flux surfaces are parameterized by a normalized flux surface coordinate $s \in [0,1]$ with $0$ being the magnetic axis and $1$ being the outermost flux surface. On each flux surface with a fixed $s$, the two-dimensional eigenvalue problem is solved.
\textblue{ In Figure~\ref{fig:varcoefs}, the variable coefficients are plotted for the flux surface $s=0.8$, together with the field-aligned mesh.}
\\
%
%In this section, we explore the numerical results of ADG for $\vB$ where metric terms from flux surfaces of MHD equilibria are included in the functions $\alpha \left( \vx \right)$ and $\beta \left( \vx \right)$. As a test case, we use the \texttt{VMEC}-equilibrium of the W7-X high-mirror case \cite[Table IV]{vmec_equilibrium}, with units normalized in Tesla and meter.\\  
\textblue{ The mapping of the flux surface between logical and physical space is shown in Figure~\ref{fig:MHD_mesh_torus}. We propose using a mesh with toroidally non-conforming interfaces and use $p_\xi = 3$ and $p_\eta  = 7$. }
\begin{figure}[!htb]
	\centering
	\subfigure[logical domain discretization]{
	\begin{tikzpicture}[line width =1pt]
		\draw (0.2, 0.5) node[inner sep=0]{
			\includegraphics[trim=0 0 0 0,clip,width=0.36\textwidth]{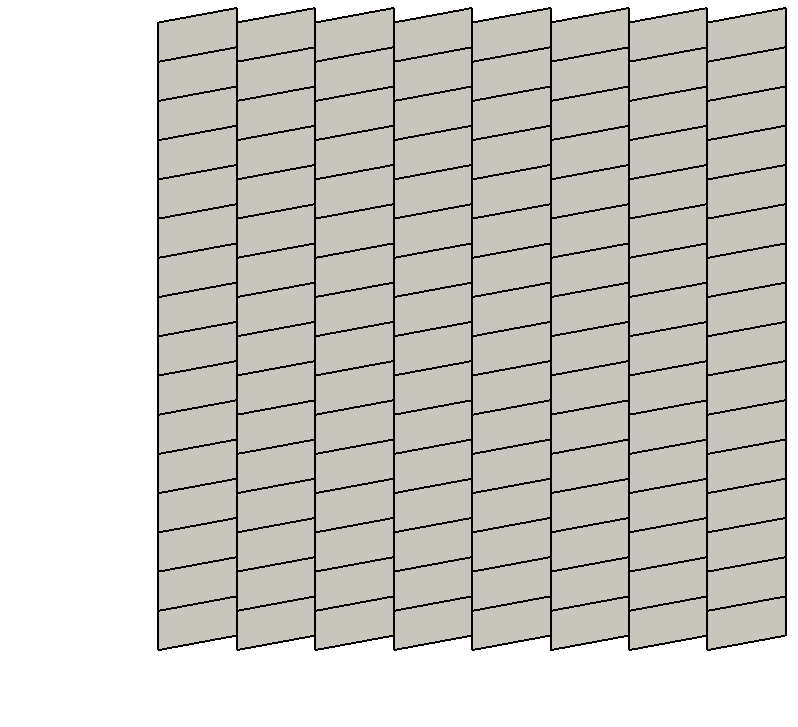}
		};
		\draw[<->] (-2,-0.7) -- (-2,1.5);
		\draw[left] (-2,0.4) node {$N_y$};
		\draw[<->] (-0.5,-2) -- (1.5,-2);
		\draw[below] (0.5,-2) node {$N_x$};
		
		\draw[->] (-2.,-2.25) -- (-2.,-1.5);
		\draw[right] (-2.,-1.5) node {$\theta$};
		\draw[->] (-2.,-2.25) -- (-1.25,-2.25);
		\draw[right] (-1.25,-2.25) node {$\varphi$};
	
	\end{tikzpicture}	
	
}
\subfigure[physical domain visualization (dark grey: one field period)]{
	\begin{tikzpicture}[line width =1pt]
	\draw (-0.2, 1.6) node[inner sep=0] {
		\includegraphics[width=0.58\textwidth]{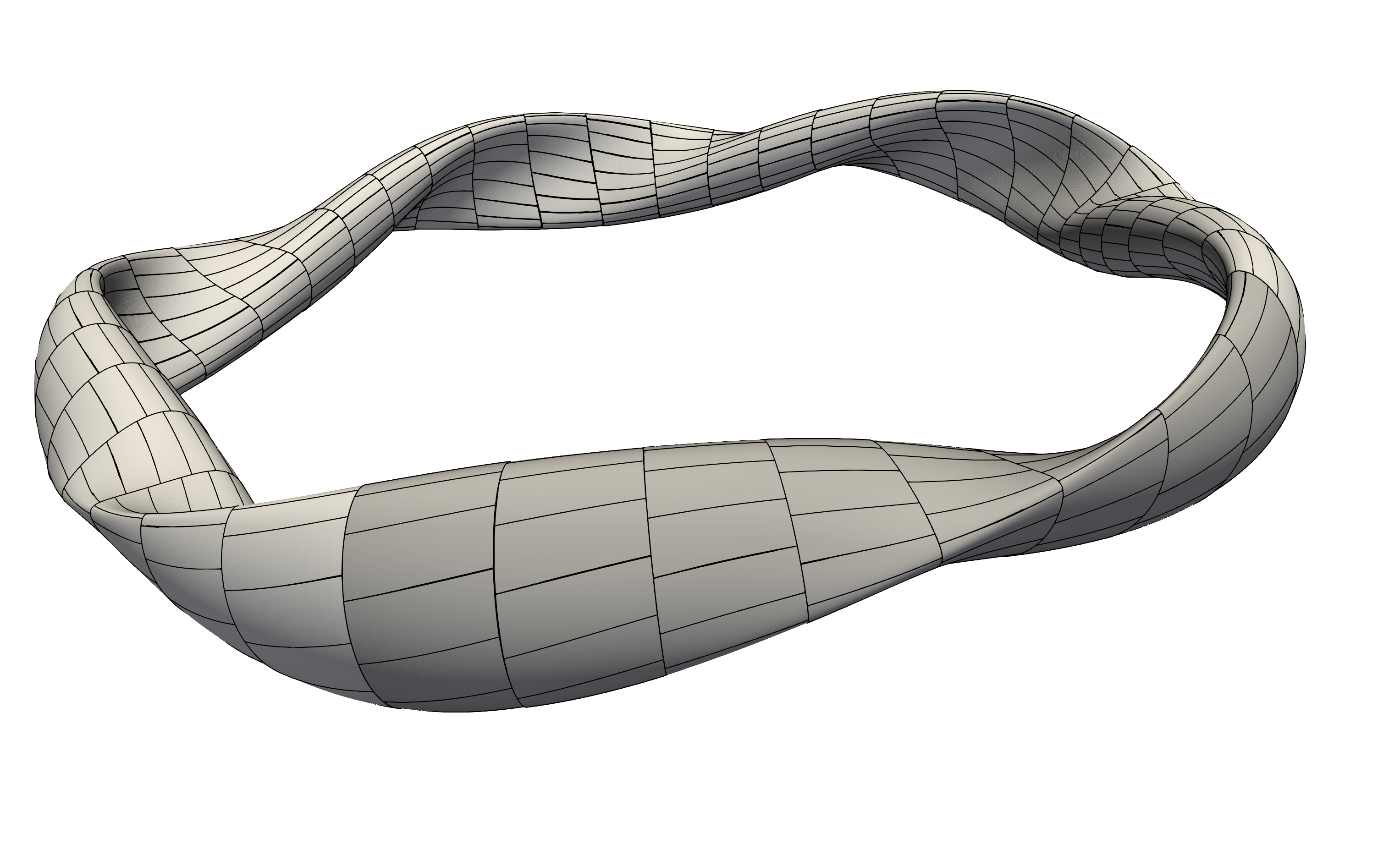}
	};
% 	\draw[->] (-0.5,1-0.866) arc (-120:0:1);	
% 	\draw[right] (1,1) node {$\varphi$};
% 	\draw[line width=5pt,draw=white,->] (0.5,0) arc (-90:180:0.33);
% 	\draw[->] (0.5,0) arc (-90:180:0.33);
% 	\draw[left] (0.15,0.5) node {$\theta$};
	\draw[->] (1.9,1-0.866) arc (-120:0:1);	
	\draw[right] (3.,1) node {$\varphi$};
	\draw[line width=5pt,draw=white,->] (2.9,0) arc (-90:180:0.33);
	\draw[->] (2.9,0) arc (-90:180:0.33);
	\draw[left] (2.6,0.5) node {$\theta$};
	
	\end{tikzpicture}
	
} 
	\caption{Locally field-aligned mesh with $N_x = 8$, $N_y = 16$ in a single field period yielding toroidally non-conforming interfaces on a flux surface for the example of a W7-X-like MHD equilibrium.}
	\label{fig:MHD_mesh_torus}
\end{figure}
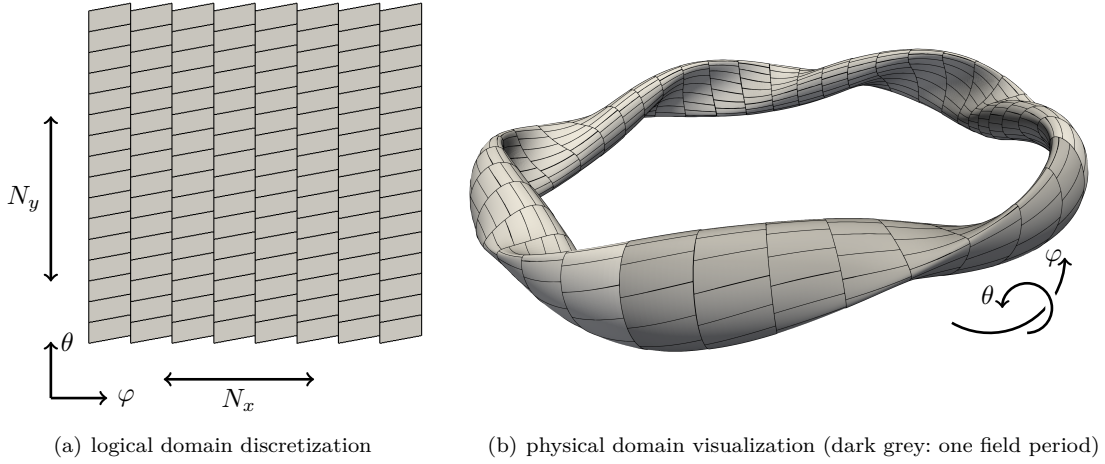\\
The number of field periods for the considered W7-X-like equilibrium is $5$. We compare results by the mesh resolution $(N_x,N_y)$ of a single field period. The total number of cells for a discretization is then $(5 N_x,N_y)$.\\
We first consider the convergence of ADG for MHD equilibria in Section~\ref{sec:ADM_convergence}. In Section~\ref{sec:ADM_alignment}, we examine the impact of the local alignment of mesh and basis. We then compare the results of ADG to an existing Fourier method in Section~\ref{sec:ADM_conti_comparison}.

\subsection{Numerically observed convergence}
\label{sec:ADM_convergence}

As no analytic results are available, we consider the discrete results of the method and their behaviour when increasing the mesh resolution.
\begin{figure}[!htb]
	\centering
	\includegraphics[width=0.98\textwidth]{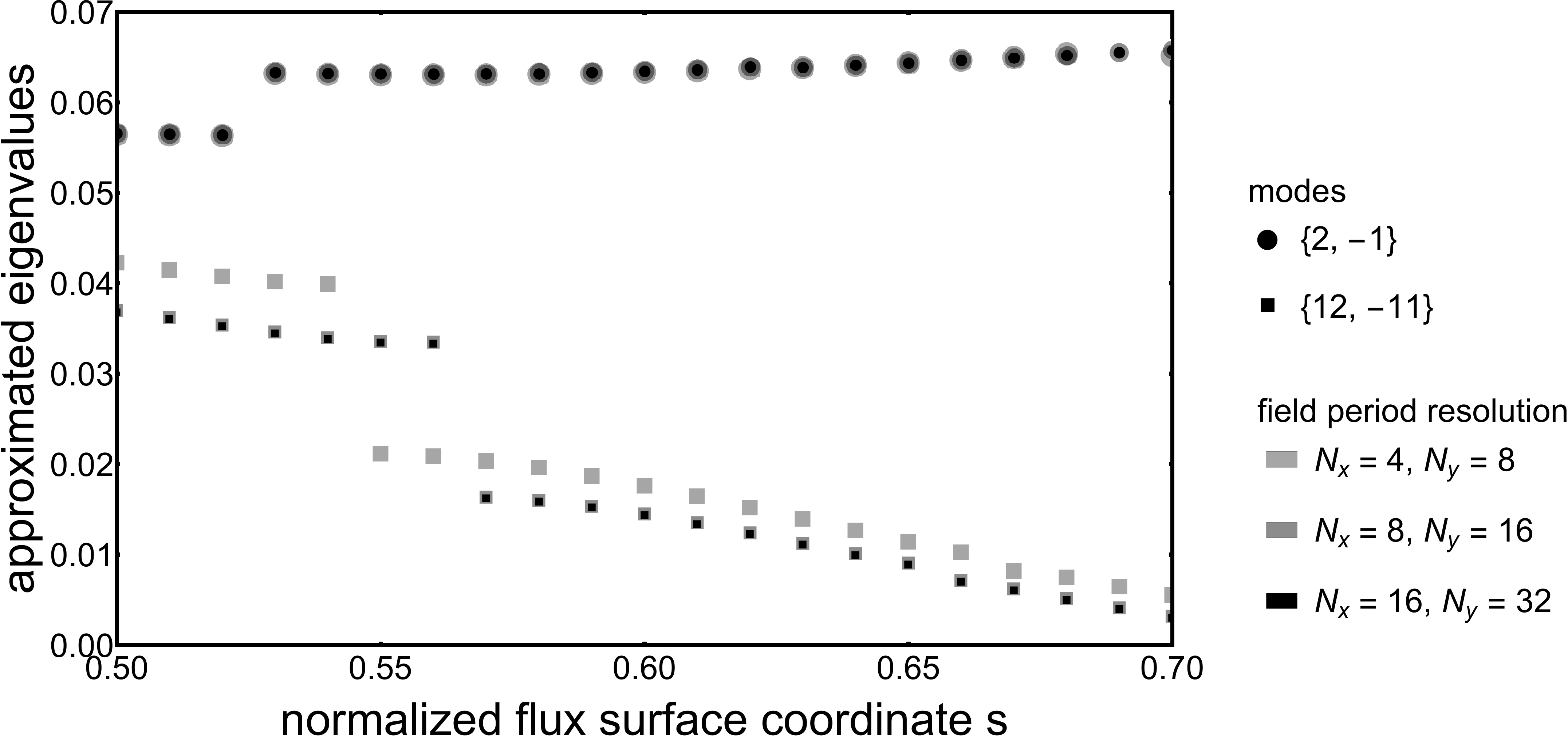}
	\caption{Convergence results of ADG with $p_\xi = 3, p_\eta = 7$ on flux surfaces of a W7-X-like MHD equilibrium. Toroidally non-conforming meshes of increasing resolution satisfying  $\dof_\perp / \dof_\parallel = 4$ within a field period are used. Modes are distinguished by shape, and resolution by colour.}	\label{fig:ADM_convergence}
\end{figure}
Figure~\ref{fig:ADM_convergence} shows the results of ADG for a selection of modes fulfilling $\dof_\perp / \dof_\parallel = 4$ within a field period. This selection is chosen to make the convergence process visible. We observe that the spectra coincide for $(2,-1)$ (round markers) for all resolutions. For $(12,-11)$, the spectra of $N_x = 8, N_y = 16$ and $N_x = 16, N_y = 32$ coincide (square markers). Therefore, we deduce that further refining the mesh yields the same results. For $N_x = 4, N_y = 8$, the eigenvalues of $(12,-11)$ are not converged yet (light gray square markers).\\
An explanation of the jumps in the spectrum is given in Section~\ref{sec:ADM_conti_comparison}.

\subsection{Impact of the local alignment}
\label{sec:ADM_alignment}

Having established the convergence of ADG, we can now compare the converged result with the non-aligned cartesian mesh. 
\begin{figure}[!htb]
	\centering
	\includegraphics[width=0.98\textwidth]{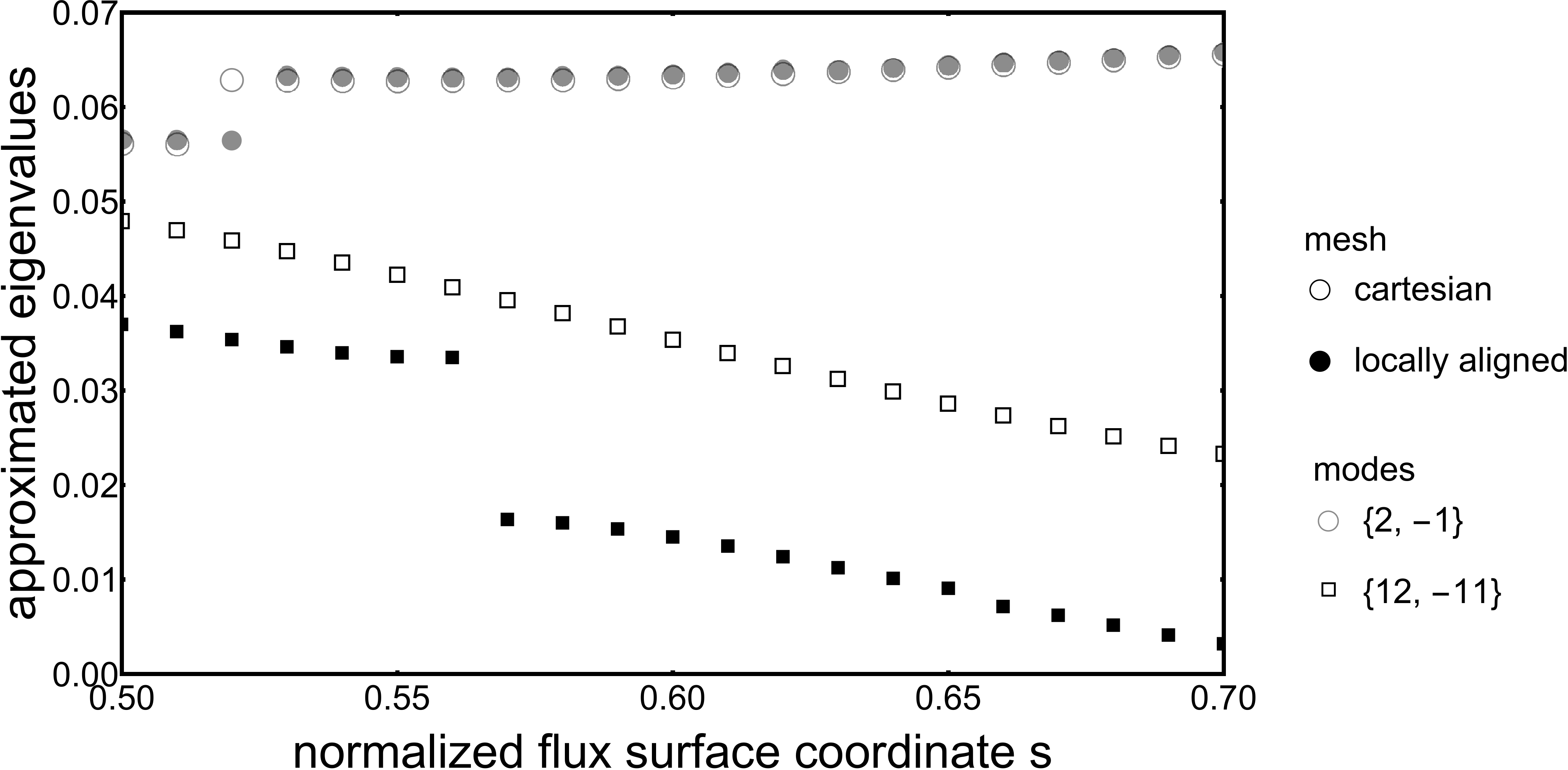}
	\caption{Results of ADG with $p_\xi = 3, p_\eta = 7, N_x = 8, N_y = 16$ using toroidally non-conforming interfaces in comparison to a non-aligned cartesian case with $p_\xi = 7, p_\eta = 7, N_x = 8, N_y = 8$ on flux surfaces of a W7-X-like MHD equilibrium. Both setups fulfill $\dof = 2^{12}$ within a field period. Modes are distinguished by shape and color, meshes are indicated by fillings.} 
	\label{fig:ADM_tor_vs_cartesian}
\end{figure}\\
Figure~\ref{fig:ADM_tor_vs_cartesian} shows these results for the same selection of modes as in Figure~\ref{fig:ADM_convergence}. The results for the low mode $(2,-1)$ coincide up to one datapoint, whereas the results for the high mode $(12,-11)$ differ by a huge margin. For the same field period resolution of $\dof = 2^{12}$, the cartesian case is not converged yet for the high mode number. 

%To reach the same level of accuracy, the cartesian mesh would need to be refined by a factor of $6$ in both directions.

%The number of degrees of freedom needed to resolve a mode number $m$ would increase quadratically in the cartesian case since both 
%
%As parallel and perpendicular resolution are not decoupled for 

%In the cartesian case, we would have to refine the mesh in both directions to resolve the high mode $(12,-11)$. %Thus, in general, the required resolution would scale quadratically with the largest mode number.

%On the other hand, the locally aligned mesh already provides a converged solution with the coarse mesh.

%We conclude that less resolution is needed when using the locally field aligned mesh with distributed resolution. 

\subsection{Comparison with a spectral code}
\label{sec:ADM_conti_comparison}

This section compares results of ADG to the spectral code CONTI \cite{CONTI} which operates on a single field period. The results of CONTI were kindly provided by Axel Könies. The setup of CONTI uses poloidal mode numbers with $m_\text{max} = 58$ and toroidal mode numbers 
\begin{equation}
n \in \left\{ -45, -40, -35, -30, \dots , 30, 35, 40, 45\right\}-\sigma 
\end{equation}
where $\sigma \in \left\{0,1, \dots, 4\right\}$ is a phase factor shift to account for different families of toroidal mode numbers. For evaluating metric terms, a resolution of $240 \times 80$ points in poloidal, toroidal direction was used. We setup ADG using a toroidally non-conforming mesh with the resolution on a single field period being $N_x = 8, N_y = 32$ and basis degrees $p_\xi = 3, p_\eta = 7$ which yields $\dof_\perp / \dof_\parallel = 8$.
\begin{figure}[!htb]
	\centering 	
	\includegraphics[width=0.98\textwidth]{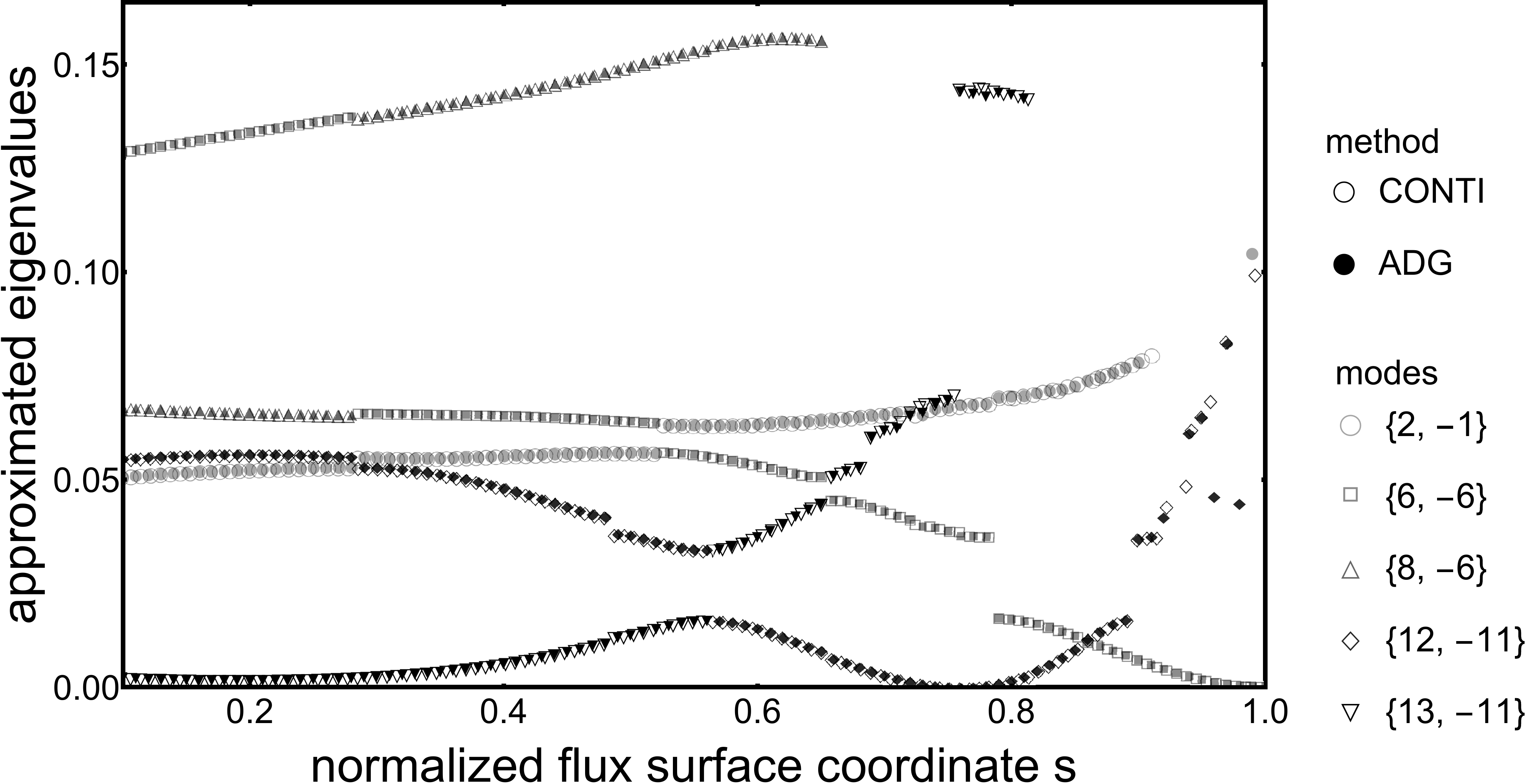}		
	\caption{Comparison of ADG using a mesh with toroidally non-conforming meshes and $p_\xi = 3, p_\eta = 7, N_x = 8, N_y = 32$ such that $\dof_\perp / \dof_\parallel = 8$ with CONTI on flux surfaces of a W7-X-like MHD equilibrium. Modes are indicated by shape and color whereas the choice of method is indicated by fillings.}
	\label{fig:ADM_Conti_comparison}
\end{figure}
Figure~\ref{fig:ADM_Conti_comparison} shows results for a larger selection of modes for ADG and CONTI. We use filled markers for ADG and empty markers for CONTI. We observe that the results of the methods overall coincide with some larger deviations for $s \geq 0.9$ at the boundary of the equilibrium where strong metric terms reside.\\
Figure~\ref{fig:ADM_Conti_comparison} allows to explain the jumps in the spectrum when tracing the eigenvalues of a single mode. As the exact eigenfunctions of the anisotropic wave equations on flux surfaces of a three-dimensional MHD equilibrium are a combination of different Fourier modes, the mode association of Section~\ref{sec:eigenvector_postprocessing} might shift when traversing from one flux surface to another, as a different Fourier mode might become the dominating part of the eigenfunction. For example, at $s\approx 0.57$ associations of $(12,-11)$ and $(13,-11)$ switch. The same holds for $(6,-6)$ and $(13,-11)$ at $s\approx 0.66$ or $(8,-6)$ and $(6,-6)$ at $s\approx 0.29$. The jumps in the spectrum are physically meaningful and are related to plasma instabilities \cite{freidberg_gaps}.

\section{Conclusion}
\label{sec:summary}

We constructed and analyzed a discontinuous Galerkin method relying on a variational mixed form of an anisotropic wave equation and operating on a non-conforming locally field-aligned mesh.\\
The numerical results of Sections~\ref{sec:results} and \ref{sec:ADM_MHD} confirm that the local alignment of mesh and basis decouples the resolution in parallel and perpendicular direction. This allows to resolve eigenfunctions with high mode numbers while providing the possibility to coarsely discretize close to constant parts. Furthermore, the size of eigenvalue errors now correlates to both the size of the mode numbers of the associated Fourier eigenmode and the size of the eigenvalue itself.\\
We examine the impact of the local alignment of mesh and basis on the spectrum. For constant coefficients, ADG yields an improvement of up to $7$ orders of magnitude in accuracy compared to a non-aligned cartesian case with the same number of degrees of freedom. A large gain in accuracy is particularly found for high mode numbers. \\
For the anisotropic wave equation with variable coefficients, modeling the flux surfaces of a three-dimensional MHD equilibrium, we  again  assert that the aligned meshes yield superior results when compared to a non-aligned cartesian mesh with the same total resolution. Furthermore, we show the convergence of ADG. Finally, the converged result of ADG for the eigenvalue spectrum shows excellent agreement with a highly resolved spectral code. 
%reproduces the same spectrum as the Fourier method CONTI. \TODO{ADG is therefore successfully benchmarked?}\\
%We conclude that the local alignment of mesh and basis allows the distribution of resolution and therefore improves the accuracy by multiple orders of magnitude in comparison to non-aligned meshes. ADG offers the flexibility to focus the computational effort on the numerically challenging and structurally important aspects of the physical problem.\\
%\TODO{Outlook on three-dimensional locally field aligned meshes for MHD equilibria?}

\subsection*{Acknowledgements}
This work has been carried out within the framework of the EUROfusion Consortium and has received funding from the Euratom research and training programme 2014-2018 under grant agreement No. 633053. The views and opinions expressed herein do not necessarily reflect those of the European Commission.\\
We want to thank Eric Sonnendrücker, Ralf Kleiber and Axel Könies for fruitful discussions and their valuable input.

\bibliographystyle{plain}

\bibliography{biblio}

\vspace*{\fill}

%\vspace*{\fill} 

% History dates

%\received{November 2012}{?}{?}

\end{document}